%% file: paper.tex
\documentclass[]{tCTM2e}
\usepackage{amsmath,amssymb}
\usepackage{booktabs}
\usepackage{color}
\usepackage{listings}

\allowdisplaybreaks[1]

\renewcommand{\vec}[1]{\mbox{\boldmath$#1$}}
\newcommand{\tensor}[1]{\mbox{\boldmath{\ensuremath{#1}}}}

\usepackage{epsfig}

\usepackage{url}

\newcommand{\sfrac}[2]{\mathchoice
  {\kern0em\raise.5ex\hbox{\the\scriptfont0 #1}\kern-.15em/
   \kern-.15em\lower.25ex\hbox{\the\scriptfont0 #2}}
  {\kern0em\raise.5ex\hbox{\the\scriptfont0 #1}\kern-.15em/
   \kern-.15em\lower.25ex\hbox{\the\scriptfont0 #2}}
  {\kern0em\raise.5ex\hbox{\the\scriptscriptfont0 #1}\kern-.2em/
   \kern-.15em\lower.25ex\hbox{\the\scriptscriptfont0 #2}}
  {#1\!/#2}}

\begin{document}

\markboth{Emmett et al.}{SMC}

\title{High-Order Algorithms for Compressible Reacting Flow with Complex Chemistry}

\author{M. Emmett$^{\rm a}$$^{\ast}$,\thanks{$^\ast$Corresponding author. Email: MWEmmett@lbl.gov\vspace{6pt}}
        W. Zhang$^{\rm a}$, and
        J. B. Bell$^{\rm a}$ \\\vspace{6pt}
$^{\rm a}${\em{Center for Computational Sciences and Engineering,
              Lawrence Berkeley National Laboratory,
              Berkeley, CA 94720.}};
}

\maketitle

\begin{abstract}
In this paper we describe a numerical algorithm for integrating the multicomponent,
reacting, compressible Navier-Stokes equations, targeted for direct numerical
simulation of combustion phenomena.  The algorithm addresses two shortcomings
of previous methods.  First, it incorporates an eighth-order
narrow stencil approximation of diffusive terms that reduces the communication compared
to existing methods and removes the need to use a filtering algorithm to remove
Nyquist frequency oscillations that are not damped with traditional approaches.  The
methodology also incorporates a multirate temporal integration strategy that provides
an efficient mechanism for treating chemical mechanisms that are stiff relative to fluid
dynamical time scales.  The overall methodology is eighth order in space with options for fourth order to eighth order in time.
The implementation uses a hybrid programming model designed for effective utilization of
many-core architectures.
We present numerical results demonstrating the convergence properties of the algorithm
with realistic chemical kinetics
and illustrating its performance characteristics.
We also present a validation example showing that the algorithm matches detailed results
obtained with an established low Mach number solver.
\bigskip

\begin{keywords}
DNS; spectral deferred corrections; finite difference; high-order methods; flame simulations
\end{keywords}\bigskip

\end{abstract}

\section{Introduction}
The simulation of a turbulent reacting flow without the use of either
turbulence models or closure models for turbulence chemistry interaction is
referred to as a direct numerical simulation, or DNS.
One of the standard computational tools used in DNS studies of combustion is
the numerical integration of the multicomponent, reacting compressible Navier-Stokes
equations.
DNS approaches have been in active use in the combustion community for more than
twenty-five years; a comprehensive review of the literature is outside the scope
of this article.  Examples of early work in the field can be found in
\cite{BaumPoinsot1994,Haworth1992,Poinsot1995,VervischPoinsot:1998,JenkinsCant:1999,ChenIm1998}.
With the advent of massively parallel, high performance computer architectures,
it is now feasible to perform fully compressible DNS simulations
for turbulent flows with detailed kinetics.  See, for example,
\cite{Tanahashi2000,Sankaran2007a,Yoo2011,Kolla2012}.

Compressible DNS approaches place a premium on high order of accuracy.
Spatial discretizations are typically based on defining high-order spatial derivative
operators, $D_x$, $D_y$ and $D_z$.
Common choices are central difference approximations \cite{Kennedy1994} or
compact differences based on Pade approximations \cite{Lele1992}.
Two codes that have seen significant application in combustion,
SENGA {\cite{JenkinsCant:1999,ChakrabortyCant:2004}
and S3D
\cite{Sankaran2007,Yoo2011,Kolla2012},
use tenth and eighth order central differences, respectively.
Second-order derivative terms are approximated using repeated application of the
first-order derivative operators, namely
\[
(a(x)U_x)_x \approx D_x ( a D_x U)
\]
This treatment of second-order derivatives does not provide adequate damping for
Nyquist-frequency oscillations. Consequently some type of higher-order filter
is employed periodically to damp high-frequency components of the system.
Alternatively, the chain rule can be applied to second-order
derivative terms.  However, the chain rule approach is known to have
stability issues (see e.g., \cite{Pantano2010}).
Discretization of the spatial operators results in a system of ordinary differential
equations (ODEs), which are typically integrated using high-order explicit Runge Kutta schemes such
as those discussed by Kennedy {\it et al.} \cite{Kennedy2001}.

Although these types of compressible DNS codes have been highly
successful, they suffer from two basic weaknesses.  First, the use of
explicit integration techniques limits the time step to the minimum of
acoustic and chemical time scales.  (Diffusive time scales are
typically less restrictive than acoustic time scales for DNS.)  This
can be particularly restrictive for detailed chemical kinetics that
can be stiff on acoustic time scales.  The time step restriction can
be relaxed by using an implicit / explicit temporal integration scheme
such as IMEX but at the expense of solving linear systems to
implicitly advance the kinetics at acoustic time scales.

The other issue with traditional compressible DNS methodology is the
evaluation of second-order derivative terms by repeated application of
the high-order first derivative operators.  This approach is efficient
in terms of floating point work but requires either multiple
communication steps or much wider layers of ghost cells for each
evaluation of the spatial operators.  Although this type of design has
been effective, as we move toward new many-core architectures that
place a premium on reduced communication, the algorithm becomes
communication bound and computational efficiency is reduced.

In this paper, we introduce two new algorithmic features that address
the weaknesses discussed above.  For the spatial discretization, we
introduce an extension of the sixth-order narrow stencil finite
difference algorithm of Kamakoti and Pantano \cite{Pantano2010} that
provides a single discretization of $(a(x)U_x)_x$ to eighth-order
accuracy.  This approach requires significantly more floating point
work than the standard approach but allows the spatial discretization
to be computed with less communication.  Furthermore, it eliminates
the need for a filter to damp Nyquist-frequency oscillations.  We also
introduce a multirate time stepping strategy based on spectral
deferred corrections (SDC) that enables the chemistry to be advanced
on a different time scale than the fluid dynamics.  This allows the
fluid dynamics to be evaluated less frequently than in a fully coupled
code resulting in shorter run times.  These algorithms are implemented
in a hybrid OpenMP/MPI parallel code named {\tt SMC}.

In \S\ref{sect:spatial} we introduce the spatial discretization used
by {\tt SMC} and the eighth-order narrow stencil discretization
of variable coefficient diffusion
operators.  In \S\ref{sect:temporal} we discuss the temporal
discretization used by {\tt SMC} and the multirate spectral deferred
correction algorithm.  In \S\ref{sect:parallel} we discuss some of the issues
associated with a hybrid parallel implementation of {\tt SMC}.  In
particular, we discuss a spatial blocking strategy that provides
effective use of OpenMP threads.  In \S\ref{sect:results} we present
numerical results illustrating the convergence properties and parallel
performance of {\tt SMC}, and show its application to dimethyl
ether-air combustion.  Finally, the main results of the paper are
summarized in \S\ref{sec:Conclusions}.

\section{Spatial Discretization}
\label{sect:spatial}

\input{smc.tex}

\section{Temporal Discretization}
\label{sect:temporal}

\input{sdc.tex}

\section{Parallelization}
\label{sect:parallel}

\input{parallel.tex}

\section{Numerical Results}
\label{sect:results}

\input{results.tex}

\section{Conclusions}\label{sec:Conclusions}

We have described a numerical algorithm for integrating the
multicomponent, reacting, compressible Navier-Stokes equations that is
eight-order in space with fourth order to eighth order temporal
integration.  The methodology uses a narrow stencil approximation of
diffusive terms that reduces the communication compared to existing
methods and removes the need to use a filtering algorithm to remove
Nyquist frequency oscillations that are not damped with traditional
approaches.  It also incorporates a multirate temporal integration
strategy that provides an efficient mechanism for treating chemical
mechanisms that are stiff relative to fluid dynamical time scales.
The narrow stencil uses more floating point operations per stencil
evaluation than the standard wide stencil but requires less
communication.  A strong scaling study demonstrated that the narrow
stencil scales better than the wide stencil which has important
ramifications for stencil performance on future multicore machines.
The multirate integration approach supports two modes for
operation. One treats the reactions explicitly using smaller time
steps than are used for the fluid mechanics; the other treats reaction
implicitly for case in which the chemical time scales are
significantly faster than the fluid dynamical time scales.  A
performance study of the multi-rate integration scheme with realistic
chemical kinetics showed that the multi-rate integrator can operate
with a larger time-step compared to single-rate integrators and
ultimately obtained an accurate solution in less time.  The
implementation uses a hybrid programming model designed for effective
utilization of a large number of threads that is targeted toward next
generation many-core architectures.  The code scales well to
approximately 100K cores and shows excellent thread performance on a
61-core Intel Xeon Phi coprocessor.  We present numerical results
demonstrating the convergence properties of the algorithm with
realistic chemical kinetics and validating the algorithm against
results obtained with an established low Mach number solver.

\section*{Acknowledgements}
This work was supported by the Applied Mathematics Program and
the Exascale Co-Design Program of the DOE Office of Advanced
Scientific Computing Research under the U.S. Department of Energy
under contract DE-AC02-05CH11231.  This research used resources of the
National Energy Research Scientific Computing Center, which is
supported by the Office of Science of the U.S. Department of Energy
under Contract No. DE-AC02-05CH11231.

\bibliographystyle{siam}
\bibliography{combustion,sdc,dns,boundary,ws}

\appendices
\section{Eighth-Order Narrow Stencil: Free Parameters and Numerical Experiments}
\label{sect:app}
\input{stencil.tex}

\end{document}

%% file: smc.tex
\subsection{Navier-Stokes Equations}
\label{sec:NSeqns}

The multicomponent reacting compressible Navier-Stokes equations are given by
\begin{subequations}
\label{eq:cns}
\begin{align}
\frac{\partial \rho}{\partial t} + \nabla \cdot (\rho
    \vec{u})= { } & 0, \\
\frac{\partial \rho \vec{u}}{\partial t} + \nabla \cdot (\rho
    \vec{u}\vec{u}) + \nabla p= { } & \nabla \cdot
  \tensor{\tau}, \\
\frac{\partial \rho Y_k}{\partial t} + \nabla \cdot (\rho Y_k \vec{u})
= { } & \rho \dot{\omega}_k  - \nabla \cdot \mathcal{F}_k, \label{eq:NS:rY} \\
\frac{\partial \rho E}{\partial t} + \nabla \cdot [(\rho E + p)
  \vec{u}] = { } & \nabla \cdot (\lambda \nabla T) + \nabla \cdot
  (\tensor{\tau} \cdot \vec{u}) - \nabla \cdot \sum_k \mathcal{F}_k h_k, \label{eq:NS:E}
\end{align}
\end{subequations}
where $\rho$ is the density, $\vec{u}$ is the velocity, $p$ is the
pressure, $E$ is the specific energy density (kinetic energy plus
internal energy), $T$ is the temperature, $\tensor{\tau}$ is the
viscous stress tensor, and $\lambda$ is the thermal conductivity.  For
each of the chemical species $k$, $Y_k$ is the mass fraction,
$\mathcal{F}_k$ is the diffusive flux, $h_k$ is the specific enthalpy, and
$\dot{\omega}_k$ is the production rate.  The system is closed by an
equation of state that specifies $p$ as a function of $\rho$, $T$ and
$Y_k$. For the examples presented here, we assume an ideal gas
mixture.

The viscous stress tensor is given by
\begin{equation}
  \tau_{ij} = \eta \left(\frac{\partial u_i}{\partial x_j} +
    \frac{\partial u_j}{\partial x_i} - \frac{2}{3}
    \delta_{ij} \nabla \cdot \vec{u} \right) +
  \xi \delta_{ij} \nabla \cdot \vec{u},
\end{equation}
where $\eta$ is the shear viscosity and $\xi$ is the bulk viscosity.

A mixture model for species diffusion is employed in {\tt SMC}.  We
define an initial approximation $\bar{\mathcal{F}}_k$ to the species
diffusion flux given by
\begin{equation}
\bar{\mathcal{F}}_k
  = - \rho D_k \left( \nabla X_k + (X_k-Y_k) \frac{\nabla p}{p} \right)
\label{eq:fbar}
\end{equation}
where $X_k$ is the mole fraction, and $D_k$ is the diffusion
coefficient of species $k$.  For species transport to be
consistent with conservation of mass we require
\begin{equation}
  \label{eq:fcons}
  \sum_k \mathcal{F}_k = 0.
\end{equation}
In general, the fluxes defined by (\ref{eq:fbar}) do not satisfy this
requirement and some type of correction is required.  When there is a
good choice for a reference species that has a high molar
concentration throughout the entire domain, one can use
(\ref{eq:fbar}) to define the species diffusion flux (i.e.,
$\mathcal{F}_k \equiv \bar{\mathcal{F}}_k $) for all species except
the reference, and then use (\ref{eq:fcons}) to define the flux of the
reference species.  This strategy is often used in combustion using
$\mathrm{N_2}$ as the reference species.  An alternative approach
\cite{Poinsot2001} is to define
a correction velocity $\vec{V}_c$ given by
\begin{equation}
  \vec{V}_c = \sum_\ell \bar{\mathcal{F}}_\ell \label{eq:Vc}
\end{equation}
and then to set
\begin{equation}
\mathcal{F}_k = \bar{\mathcal{F}}_k - Y_k \vec{V}_c.
\end{equation}
The {\tt SMC} code supports both of these operating modes.  We will
discuss these approaches further in \S\ref{sec:correction}.

\subsection{Eighth-Order and Sixth-Order Discretization}
\label{sect:stencil}

The {\tt SMC} code discretizes the governing Navier-Stokes equations
(\ref{eq:cns}) in space using high-order centered finite-difference
methods.  For first order derivatives a standard $8^{\rm th}$ order
stencil is employed.  For second order derivatives with variable
coefficients of the form
\begin{equation}
  \label{eq:vc2}
  \frac{\partial}{\partial x} \left(a(x) \frac{\partial u}{\partial x}
\right),
\end{equation}
a novel $8^{\rm th}$ order narrow stencil is employed.  This stencil
is an extension of those developed by Kamakoti and Pantano
\cite{Pantano2010}.  The derivatives in (\ref{eq:vc2}) on a uniform
grid are approximated by
\begin{equation}
\frac{\partial}{\partial x} \left(a \frac{\partial u}{\partial x}
\right)\bigg{|}_{i} \approx \frac{H_{i+1/2} - H_{i-1/2}}{\Delta x^2}, \label{eq:DaDu}
\end{equation}
where
\begin{equation}
  H_{i+1/2} = \sum_{m=-s+1}^{s} \sum_{n=-s+1}^{s} a_{i+m} M_{mn} u_{i+n}, \label{eq:H}
\end{equation}
$s=4$, and $M$ is an $8\times 8$ matrix (in general, for a discretization of
order $2s$, $M$ is a $2s \times 2s$ matrix).  This discretization is
conservative because (\ref{eq:DaDu}) is in flux form.  Note, however,
that $H_{i+1/2}$ is not a high-order approximation of the flux at
$x_{i+1/2}$ even though the flux difference (\ref{eq:DaDu})
approximates the derivative of the flux to order $2s$.  This type of
stencil is {\emph{narrow}} because the width of the stencil is only
$2s+1$ for both $u$ and $a$.  In contrast, the width of a wide
stencil in which the first-order derivative operator is applied twice
is $4s+1$ for $u$ and $2s+1$ for $a$.

Following the strategy in \cite{Pantano2010}, a family of $8^{\rm th}$
order narrow stencils for the second order derivative in
(\ref{eq:vc2}) has been derived.  The resulting $8^{\rm th}$ order stencil
matrix is given by
\begin{equation}
  M = \left[ \begin{array}{cccccccc}
      m_{11} & m_{12} & m_{13} & m_{14} & m_{15} & 0 & 0 & 0\\
      m_{21} & m_{22} & m_{23} & m_{24} & m_{25} & m_{26} & 0 & 0 \\
      m_{31} & m_{32} & m_{33} & m_{34} & m_{35} & m_{36} & m_{37} & 0 \\
      m_{41} & m_{42} & m_{43} & m_{44} & m_{45} & m_{46} & m_{47} & m_{48} \\
      -m_{48} & -m_{47} & -m_{46} & -m_{45} & -m_{44} & -m_{43} & -m_{42} & -m_{41} \\
      0 & -m_{37} & -m_{36} & -m_{35} & -m_{34} & -m_{33} & -m_{32} & -m_{31} \\
      0 & 0 & -m_{26} & -m_{25} & -m_{24} & -m_{23} & -m_{22} & -m_{21} \\
      0 & 0 & 0 & -m_{15} & -m_{14} & -m_{13} & -m_{12} & -m_{11}
      \end{array} \right],
\end{equation}
where
\begin{align*}
  m_{11} &= \frac{5}{336} + m_{48}, &
  m_{12} &= -\frac{83}{3600} - \frac{1}{5} m_{47} - \frac{14}{5} m_{48},\\
  m_{13} &= \frac{299}{50400} + \frac{2}{5} m_{47} + \frac{13}{5} m_{48}, &
  m_{14} &= \frac{17}{12600} - \frac{1}{5} m_{47} - \frac{4}{5} m_{48},\\
  m_{15} &= \frac{1}{1120}, &
  m_{21} &= -\frac{11}{560} - 2 m_{48},\\
  m_{22} &= -\frac{31}{360} + m_{47} + 3 m_{48}, &
  m_{23} &= \frac{41}{200} - \frac{9}{5} m_{47} + \frac{4}{5} m_{48},\\
  m_{24} &= -\frac{5927}{50400} + \frac{4}{5} m_{47} - \frac{9}{5} m_{48}, &
  m_{25} &= \frac{17}{600} - \frac{1}{5} m_{47} - \frac{4}{5} m_{48},\\
  m_{26} &= -\frac{503}{50400} + \frac{1}{5} m_{47} + \frac{4}{5} m_{48}, &
  m_{31} &= -\frac{1}{280},\\
  m_{32} &= \frac{1097}{5040} - 2 m_{47} + 6 m_{48}, &
  m_{33} &= -\frac{1349}{10080} + 3 m_{47} - 12 m_{48},\\
  m_{34} &= -\frac{887}{5040} - m_{47} + 6 m_{48}, &
  m_{35} &= \frac{3613}{50400} +\frac{4}{5} m_{47} - \frac{9}{5} m_{48},\\
  m_{36} &= \frac{467}{25200} - \frac{3}{5} m_{47} + \frac{18}{5} m_{48}, &
  m_{37} &= \frac{139}{25200} - \frac{1}{5} m_{47} - \frac{9}{5} m_{48},\\
  m_{41} &= \frac{17}{1680} + 2 m_{48}, &
  m_{42} &= -\frac{319}{2520} + 2 m_{47} - 8 m_{48},\\
  m_{43} &= -\frac{919}{5040} - 2 m_{47} + 6 m_{48}, &
  m_{44} &= -\frac{445}{2016},\\
  m_{45} &= \frac{583}{720} - m_{47} + 6 m_{48}, &
  m_{46} &= -\frac{65}{224} - 7 m_{48},
\end{align*}
and $m_{47}$ and $m_{48}$ are two free parameters.  Note that the
well-posedness requirement established in \cite{Pantano2010} is
satisfied by the stencil matrix above independent of the two free
parameters.  In {\tt SMC} the default values of the free parameters
were set to $m_{47} = 3557/44100$ and $m_{48} = -2083/117600$ to
minimize an upper bound on truncation error that is calculated by
summing the absolute value of each term in the truncation error.  Additional
analysis of truncation errors and numerical experiments on the $8^{\rm
  th}$ order narrow stencil using various choices of the two
parameters are presented in Appendix~\ref{sect:app}.

The {\tt SMC} code uses a $6^{\rm th}$ order stencil near physical
boundaries where the conditions outside the domain are not well
specified (see \S\ref{sec:bndry}).  The general $6^{\rm th}$ order
stencil matrix is given by
\begin{equation}
  M = \left[ \begin{array}{cccccc}
      m_{11} & m_{12} & m_{13} & m_{14} & 0     & 0 \\
      m_{21} & m_{22} & m_{23} & m_{24} & m_{25} & 0 \\
      m_{31} & m_{32} & m_{33} & m_{34} & m_{35} & m_{36} \\
      -m_{36} & -m_{35} & -m_{34} & -m_{33} & -m_{32} & -m_{31} \\
      0      & -m_{25} & -m_{24} & -m_{23} & -m_{22} & -m_{21} \\
      0      &  0     & -m_{14} & -m_{13} & -m_{12} & -m_{11}
      \end{array} \right],
\end{equation}
where
\begin{align*}
  m_{11} &= -\frac{11}{180} + m_{36}, &
  m_{12} &= \frac{1}{9} - 2 m_{36},\\
  m_{13} &= -\frac{1}{18} + m_{36}, &
  m_{14} &= \frac{1}{180},\\
  m_{21} &= \frac{7}{60} - 3 m_{36}, &
  m_{22} &= -\frac{1}{120} + 5 m_{36},\\
  m_{23} &= -\frac{17}{90} - 2 m_{36}, &
  m_{24} &= \frac{5}{72} + m_{36},\\
  m_{25} &= \frac{1}{90} - m_{36}, &
  m_{31} &= -\frac{1}{15} + 3 m_{36},\\
  m_{32} &= -\frac{11}{60} - 3 m_{36}, &
  m_{33} &= -\frac{101}{360},\\
  m_{34} &= \frac{137}{180} - 2 m_{36}, &
  m_{35} &= -\frac{83}{360} + m_{36},
\end{align*}
and $m_{36}$ is a free parameter.  The $6^{\rm th}$ order narrow
stencil developed in \cite{Pantano2010} is a special case of the
general $6^{\rm th}$ order stencil when $m_{36} = 0.220063$.
Like the $8^{\rm th}$ order
stencil, the well-posedness and the order of accuracy for the $6^{\rm
  th}$ order stencil is independent of the free parameter.  In {\tt SMC}, the free
parameter is set to $281/3600$ to minimize an upper bound
on truncation error that is calculated by summing the absolute value
of each term in the truncation error.

\subsection{Cost Analysis of Eighth-Order Stencil}
\label{sec:cost}

The $8^{\rm th}$ order narrow stencil presented in
\S~\ref{sect:stencil} requires significantly more floating-point
operations (FLOPs) than the standard wide stencil, but has the benefit
of less communication.  The computational cost of a single evaluation
the narrow stencil, for $N$ cells in one dimension, is $112 N + 111$
FLOPs (this is the cost of operations in \eqref{eq:H} and
\eqref{eq:DaDu}, excluding the division by $\Delta x^2$).  For the
$8^{\rm th}$ order wide stencil with two steps of communication, the
computational cost on $N$ cells is $23 N + 8$ FLOPs.  Hence, the
narrow stencil costs $89 N + 103$ more FLOPs than the wide stencil on
$N$ cells in one dimension.  As for communication, the wide stencil
requires $64$ more bytes than the narrow stencil for $N$ cells with 8
ghost cells in one dimension.  Here we assume that an 8-byte
double-precision floating-point format is used.  Ignoring latency, the
difference in time for the two algorithms is then
\[
\frac{89N + 103}{F} - \frac{64}{B}
\]
where $F$ is the floating-point operations per second (FLOPS) and $B$
is the network bandwidth.  To estimate this difference on current
machines, we used the National Energy Research Scientific Computing
Center's newest supercomputer Edison as a reference.  Each compute
node on Edison has a peak performance of 460.8 GFLOPs (Giga-FLOPs per
second), and the MPI bandwidth is about 8 GB/s \footnote{The
  configuration of Edison is available at
  \url{http://www.nersc.gov/users/computational-systems/edison/configuration/}.}.
For $N = 64$ in one dimension, the run time cost for the extra FLOPs
of the narrow stencil is about $1.3 \times 10^{-8}\,\mathrm{s}$, which
is about 1.6 times of the extra communication cost of the wide
stencil, $8.0 \times 10^{-9}\,\mathrm{s}$.  It should be emphasized
that the estimates here are crude for several reasons.  We neglected
the efficiency of both floating-point and network operations for
simplicity.  In counting FLOPs, we did not consider a fused
multiply-add instruction or memory bandwidth.  In estimating
communication costs, we did not consider network latency, the cost of
data movement if MPI messages are aggregated to reduce latency, or
communication hiding through overlapping communication and
computation.  Moreover, the efficiency of network operations is likely
to drop when the number of nodes used increases, whereas the
efficiency of FLOPs alone is independent of the number of total nodes.
Nevertheless, as multicore systems continue to evolve, the ratio of
FLOPS to network bandwidth is most likely to increase.  Thus, we
expect that the narrow stencil will cost less than the wide stencil on
future systems, although any net benefit will depend on the systems
configuration.  As a data point, we note that both stencils are
implemented in {\tt SMC}, and for a test with $512^3$ cells using 512
MPI ranks and 12 threads on 6144 cores of Edison, a run using the
narrow stencil took about the same time as using the wide stencil.
Further measurements are presented in \S\ref{sect:performance}.

\subsection{Correction for Species Diffusion}
\label{sec:correction}

In \S\ref{sec:NSeqns}, we introduced the two approaches employed in
{\tt SMC} to enforce the consistency of the mixture model of species
diffusion with respect to mass conservation.  The first approach uses
a reference species and is straightforward to implement.  We note that
the last term in (\ref{eq:NS:E}) can be treated as $- \nabla \cdot
\sum_k \mathcal{F}_k (h_k - h_{\mathrm{ref}})$, where
$h_{\mathrm{ref}}$ is the species enthalpy for the reference species,
because of (\ref{eq:fcons}).  The second approach uses a correction
velocity (\ref{eq:Vc}) and involves the computation of $\nabla \cdot
(Y_k \vec{V}_c)$ and $\nabla \cdot (h_k Y_k \vec{V}_c)$.  Although one
might attempt to compute the correction velocity according to
$V_{c,i+1/2} = \sum_\ell H_{\ell,i+1/2}$, where $H_{\ell,i+1/2}$ for
species $\ell$ is computed using (\ref{eq:H}), this would reduce the
accuracy for correction velocity to second-order since $H_{i+1/2}$ is
only a second-order approximation to the flux at $i+1/2$.  To maintain
$8^{\rm th}$ order accuracy without incurring extra communication and
too much extra computation, we employ the chain rule for the
correction terms.  For example, we compute $\nabla \cdot (h_k Y_k
\vec{V}_c)$ as $h_k Y_k \nabla \cdot \vec{V}_c + \nabla (h_kY_k) \cdot
\vec{V}_c$, and use $\nabla \cdot \vec{V}_c = \sum_\ell \nabla \cdot
\bar{\mathcal{F}}_\ell$, where $\nabla \cdot \bar{\mathcal{F}}_\ell$
is computed using the $8^{\rm th}$ order narrow stencil for second
order derivatives.  The velocity correction form arises naturally in
multicomponent transport theory; namely, it represents the first term
in an expansion of the full multicomponent diffusion matrix. (See
Giovangigli \cite{Giovangigli1999}.)  For that reason, velocity
correction is the preferred option in {\tt SMC}.

\subsection{Physical Boundaries}
\label{sec:bndry}

The {\tt SMC} code uses characteristic boundary conditions to treat
physical boundaries.  Previous studies in the literature have shown
that characteristic boundary conditions can successfully lower
acoustic reflections at open boundaries without suffering from
numerical instabilities.  Thompson \cite{Thompson1987,Thompson1990}
applied the one-dimensional approximation of the characteristic
boundary conditions to the hyperbolic Euler equations, and the
one-dimensional formulation was extended by Poinsot and Lele
\cite{PoinsotLele1992} to the viscous Navier-Stokes equations.  For
multicomponent reacting flows, Sutherland and Kennedy
\cite{Sutherland2003} further improved the treatment of boundary
conditions by the inclusion of reactive source terms in the
formulation.  In {\tt SMC}, we have adopted the formulation of
characteristic boundary conditions proposed in \cite{Yoo2005,Yoo2007},
in which multi-dimensional effects are also included.

The $8^{\rm{ th}}$ order stencils used in {\tt SMC} require four cells
on each side of the point where the derivative
is being calculated.  Near the physical boundaries where the conditions
outside the domain are unknown, lower order stencils are used for
derivatives with respect to the normal direction.
At the fourth and third cells from the boundary we use $6^{\rm{th}}$ and $4^{\rm{th}}$ order centered stencils
respectively, for both
advection and diffusion.
At cells next to the boundary we reduce to centered second-order for diffusion and a biased third-order
discretization for advection.
Finally, for the boundary
cells, a completely biased second-order stencil is used for
diffusion and a completely biased third-order stencil is used for advection.

%% file: sdc.tex
\subsection{Spectral Deferred Corrections}
\label{sect:sdc}

\newcommand{\R}{{\mathbb{R}}}
\newcommand{\C}{{\mathbb{C}}}
\newcommand{\Qmat}{{\bm{Q}}}
\newcommand{\Smat}{{\bm{S}}}
\newcommand{\Uvec}{{\bm{U}}}
\newcommand{\Fvec}{{\bm{F}}}
\newcommand{\tc}{{\tau}}

SDC methods for ODEs were first introduced in \cite{dutt2000spectral} and
have been subsequently refined and extended in
\cite{minion2003semi,minion2004semi,huang2006accelerating,hansen2006convergence,bourlioux2003high,layton2004conservative}.
SDC methods construct high-order solutions within one time step by
iteratively approximating a series of correction equations at
collocation nodes using low-order sub-stepping methods.

Consider the generic ODE initial value problem
\begin{equation}
  \label{eq:ode}
  u'(t) = f\bigl(u(t), t\bigr), \qquad u(0) = u_0
\end{equation}
where $t \in [0,T]$; $u_0,\, u(t) \in \C^N$; and $f: \R \times \C^N
\rightarrow \C^N$.  SDC methods are derived by considering the
equivalent Picard integral form of \eqref{eq:ode} given by
\begin{equation}
  \label{eq:picard}
  u(t) = u_0 + \int_{0}^t f\bigl(u(s), s\bigr) \,ds.
\end{equation}
A single time step $[t_n, t_{n+1}]$ is divided into a set of
intermediate sub-steps by defining $M+1$ collocation points $\tc_m \in
[t_n, t_{n+1}]$ such that $t_n = \tc_0 < \tc_1 < \cdots < \tc_M = t_{n+1}$.
Then, the integrals of $f\bigl(u(t), t\bigr)$ over each of the
intervals $[t_n, \tc_m]$ are approximated by
\begin{equation}
  \label{eq:gaussapp}
  I_n^m
    \equiv \int_{t_n}^{\tc_{m}}  f\bigl(u(s), s\bigr) \,ds
    \approx \Delta t \sum_{j=0}^M q_{mj} f(U_j, \tc_j)
\end{equation}
where $U_j \approx u(\tc_j)$, $\Delta t = t_{n+1} - t_n$, and $q_{mj}$
are quadrature weights.  The quadrature weights that give the highest
order of accuracy given the collocation points $\tc_m$ are obtained by
computing exact integrals of the Lagrange interpolating polynomial
over the collocation points $\tc_m$.  Note that each of the integral
approximations $I_n^m$ (that represent the integral of $f$ from $t_n$
to $\tc_m$) in \eqref{eq:gaussapp} depends on the function values
$f(U_m, \tc_m)$ at all of the collocation nodes $\tc_m$.

To simplify notation, we define the \emph{integration matrix} $\Qmat$
to be the $M \times (M+1)$ matrix consisting of entries $q_{mj}$; and
the vectors $\Uvec \equiv [U_1, \cdots, U_M]$ and $\Fvec \equiv [
f(U_0, t_0), \cdots, f(U_M, t_M)]$.  With these definitions, the
Picard equation \eqref{eq:picard} within the time step $[t_n,
t_{n+1}]$ is approximated by
\begin{equation}
  \label{eq:compact}
  \Uvec = \Uvec_0 + \Delta t\, \Qmat\, \Fvec
\end{equation}
where $\Uvec_0 = U_0 \bm{1}$.  Note again that the integration matrix
$\Qmat$ is dense so that each entry of $\Uvec$ depends on all other
entries of $\Uvec$ (through $\Fvec$) and $U_0$. Thus,
\eqref{eq:compact} is an implicit equation for the unknowns in
$\Uvec$ at all of the quadrature points.  Finally, we note that the solution of \eqref{eq:compact}
corresponds to the collocation or implicit Runge-Kutta solution of
\eqref{eq:ode} over the nodes $\tc_m$. Hence SDC can be considered as an
iterative method for solving the spectral collocation formulation.

The SDC scheme used here begins by spreading the initial condition
$u_n$ to each of the collocation nodes so that the provisional
solution $\Uvec^0$ is given by $\Uvec^0 = [U_0, \cdots, U_0]$.
Subsequent iterations (denoted by $k$ superscripts) proceed by
applying the \emph{node-to-node integration matrix} $\Smat$ to
$\Fvec^k$ and correcting the result using a forward-Euler time-stepper
between the collocation nodes.  The node-to-node integration matrix
$\Smat$ is used to approximate the integrals $I_m^{m+1} =
\int_{\tc_m}^{\tc_{m+1}} f(u(s), s) \,ds$ (as opposed to the
integration matrix $\Qmat$ which approximates $I_n^m$) and is
constructed in a manner similar to the integration matrix
$\Qmat$.  The update equation corresponding to the forward-Euler
sub-stepping method for computing $\Uvec^{k+1}$ is given by
\begin{equation}
  \label{eq:expsdc}
  U^{k+1}_{m+1} = U^{k+1}_m +
    \Delta t_m
      \bigl[ f(U^{k+1}_{m}, t_{m}) - f(U^k_{m}, t_{m}) \bigr]
    + \Delta t\, S^{k,m}
\end{equation}
where $S^{k,m}$ is the $m^{\rm th}$ row of $\Smat \Fvec^k$.  The
process of computing \eqref{eq:expsdc} at all of the collocation nodes
$\tc_m$ is referred to as an \emph{SDC sweep} or an \emph{SDC
  iteration}.  The accuracy of the solution generated after $k$ SDC
iterations done with such a first-order method is formally $O(\Delta
t^k)$ as long as the spectral integration rule (which is determined by
the choice of collocation nodes $\tc_m$) is at least order $k$.

Here, a method of
lines discretization based on the spatial discretization presented in
\S\ref{sect:spatial} is used to reduce the partial differential
equations (PDEs) in question to a large system of ODEs.

\newcommand{\Nfeval}{N_{\rm evals}}
\newcommand{\Niters}{K}

The computational cost of the SDC scheme used here is determined by
the number of nodes $M+1$ chosen and the number of SDC iterations
$\Niters$ taken per time step.  The number of resulting function
evaluations $\Nfeval$ is $\Nfeval = \Niters M$.  For example, with 3
Gauss-Lobatto ($M=2$) nodes, the resulting integration rule is $4^{\rm
  th}$ order.  To achieve this formal order of accuracy we perform
$\Niters=4$ SDC iterations per time step and hence we require $\Nfeval
= 8$ function evaluations (note that $F(U_0, 0)$ is recycled from the
previous time step, and hence the first time-step requires 9 function
evaluations).

By increasing the number of SDC nodes and adjusting the number of SDC
iterations taken, we can construct schemes of arbitrary order.  For
example, to construct a $6^{\rm th}$ order scheme one could use a 5
point Gauss-Lobatto integration rule (which is $8^{\rm th}$ order
accurate), but only take 6 SDC iterations.  For high order explicit
schemes used here we have observed that the SDC iterations sometimes
converge to the collocation solution in fewer iterations than is
formally required.  This can be detected by computing the SDC residual
\begin{equation}
  \label{eq:residual}
  R^k = U_0 + \bm{q} \cdot \Fvec^k - U_M^k,
\end{equation}
where $\bm{q}$ is the last row of $\Qmat$, and comparing successive
residuals.  If $|R^{k-1}| / |R^{k}|$ is close to one, the SDC
iterations have converged to the collocation solution and subsequent
iterations can be skipped.

\subsection{Multirate Integration}
\label{sect:mrsdc}

\newcommand{\fslow}{f_{1}}
\newcommand{\ffast}{f_{2}}
\newcommand{\tvec}{\bm{t}}

Multirate SDC (MRSDC) methods use a hierarchy of SDC schemes to
integrate systems that contain processes with disparate time-scales
more efficiently than fully coupled methods
\cite{bourlioux2003high,layton2004conservative}.  A traditional MRSDC
method integrates processes with long time-scales with fewer SDC nodes
than those with short time-scales.  They are similar in spirit to
sub-cycling schemes in which processes with short characteristic
time-scales are integrated with smaller time steps than those with
long characteristic time-scales.  MRSDC schemes, however, update the
short time-scale components of the solutions at the long time-scale
time steps during each MRSDC iteration, as opposed to only once for
typical sub-cycling schemes, which results in tighter coupling between
processes throughout each time-step
\cite{bourlioux2003high,layton2004conservative}.  The MRSDC schemes
presented here generalize the structure and coupling of the quadrature
node hierarchy relative to the schemes first presented in
\cite{bourlioux2003high,layton2004conservative}.

Before presenting MRSDC in detail, we note that the primary advantage
of MRSDC schemes lies in how the processes are coupled.  For systems
in which the short time-scale process is moderately stiff, using more
SDC nodes than the long time-scale process helps alleviate the
stiffness of the short time-scale processes. However, extremely
stiff processes are best integrated between SDC nodes with, for
example, a variable order/step-size backward differentiation formula
(BDF) scheme.  In this case, the short time-scale processes can in
fact be treated with fewer SDC nodes than long time-scale processes.

Consider the generic ODE initial value problem with two disparate
time-scales
\begin{equation}
  \label{eq:mrode}
  u'(t) = \fslow \bigl(u(t), t\bigr) + \ffast \bigl(u(t), t\bigr), \qquad u(0) = u_0
\end{equation}
where we have split the generic right-hand-side $f$ in~\eqref{eq:ode}
into two components: $\fslow$ and $\ffast$.  A single time step $[t_n,
t_{n+1}]$ is first divided into a set of intermediate sub-steps by
defining $M_1+1$ collocation points such that $\tc_m \in [t_n,
t_{n+1}]$.  The vector of these $M_1+1$ collocation points $\tc_m$ is
denoted $\tvec_1$, and these collocation points will be used to
integrate the $\fslow$ component.  The time step is further divided by
defining $M_2+1$ collocation points $\tc_p \in [t_n, t_{n+1}]$ such
that $\tvec_1 \subset \tvec_2$, where $\tvec_2$ is the vector of these
new $\tc_p$ collocation points.  These collocation points will be used
to integrate the $\ffast$ component.  That is, a hierarchy of
collocation points is created with $\tvec_1 \subset \cdots \subset
\tvec_N$.  The splitting of the right-hand-side in \eqref{eq:mrode}
should be informed by the numerics at hand.  Traditionally, slow
time-scale processes are split into $\fslow$ and placed on the
$\tvec_1$ nodes, and fast time-scale processes are split into $\ffast$
and places on the $\tvec_2$ nodes.  However, if the fast time-scale
processes in the MRSDC update equation \eqref{eq:mrupdate} are
extremely stiff then it is more computationally effective to advance
these processes using stiff integration algorithms over a longer time
scale.  In this case, we can swap the role of of the processes within
MRSDC and treat the explicit update on the faster time scale (i.e.,
with the $\tvec_2$ nodes).  As such, we hereafter refer to the
collocation points in $\tvec_1$ as the ``coarse nodes'' and those in
$\tvec_2$ as the ``fine nodes''.

\begin{figure}
  \centering
  \includegraphics{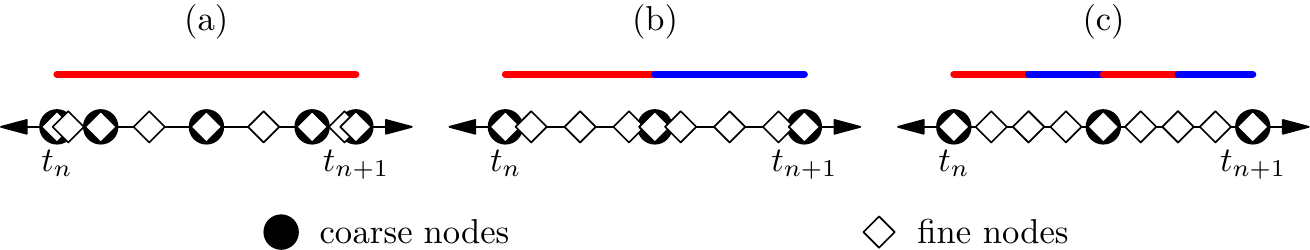}
  \caption{MRSDC quadrature nodes.  (a) The left axis shows an MRSDC
    hierarchy with $\tvec_1$ corresponding to 5 Clenshaw-Curtis nodes
    and $\tvec_2$ corresponding to 9 Clenshaw-Curtis nodes.  (b) The
    middle axis shows an MRSDC hierarchy with $\tvec_1$ corresponding
    to 3 Gauss-Lobatto nodes and $\tvec_2$ corresponding to 5
    Gauss-Lobatto nodes between each pair of nodes in $\tvec_1$.  (c)
    The right axis shows an MRSDC hierarchy with $\tvec_1$
    corresponding to 3 Gauss-Lobatto nodes and $\tvec_2$ corresponding
    to 3 Gauss-Lobatto nodes repeated twice between each pair of nodes
    in $\tvec_1$.  The red and blue lines show that the MRSDC
    integration matrix $\Qmat_{2,2}$ is constructed with (a) one
    polynomial spanning all 9 fine nodes, (b) two polynomials spanning
    5 fine nodes each, and (c) four polynomials spanning 3 fine nodes
    each.}
  \label{fig:mrpts}
\end{figure}

The SDC nodes $\tvec_1$ for the coarse component are typically chosen
to correspond to a formal quadrature rule like the Gauss-Lobatto rule.
These nodes determine the overall order of accuracy of the resulting
MRSDC scheme.  Subsequently, there is some flexibility in how the
nodes for fine components are chosen.  For example, the points in
$\tvec_2$ can be chosen to correspond to a formal quadrature rule
\begin{enumerate}
\item[(a)] over the entire time step $[t_n, t_{n+1}]$;
\item[(b)] over each of the intervals defined by the points in $\tvec_1$;
\item[(c)] over several sub-intervals of the intervals defined by the
  points in $\tvec_1$;
\end{enumerate}
Note that in the first case, a quadrature rule that nests (or embeds)
itself naturally is required (such as the Clenshaw-Curtis rule),
whereas in the second and third cases the nodes in $\tvec_2$ usually
do not correspond to a formal quadrature rule over the entire time
step $[t_n, t_{n+1}]$.  In the second and third cases we note that the
basis used to construct the $\Qmat$ and $\Smat$ integration matrices
are typically piecewise polynomial -- this allows us to shrink the
effective time-step size of the fine process without constructing an
excessively high-order polynomial interpolant.  These three
possibilities are depicted in Fig.~\ref{fig:mrpts}.

Once the SDC nodes for the desired MRSDC scheme have been chosen, the
node-to-node and regular integration matrices ($\Smat$ and $\Qmat$,
respectively) can be constructed.  While single-rate SDC schemes, as
presented in \S\ref{sect:sdc}, require one set of $\Smat$ and
$\Qmat$ matrices, the two-component MRSDC scheme used here requires
two sets of $\Smat$ and $\Qmat$ matrices.  These are denoted with
subscripts $i, j$ so that the action of $\Smat_{i,j}$ approximates the
node-to-node integrals between nodes in $\tvec_i$ of the $\Fvec_j$
component that ``lives'' on the $\tvec_j$ nodes, and similarly for the
$\Qmat_{i,j}$ matrices.  Note that both $\Smat_{i,j}$ and
$\Qmat_{i,j}$ are of size $M_i \times (M_j + 1)$.

With these definitions, the Picard equation~\eqref{eq:picard} within
the time step $[t_n, t_{n+1}]$ is approximated by
\begin{equation}
  \Uvec = \Uvec_0 + \Delta t\, \Qmat_{2,1} \Fvec_1 + \Delta t\, \Qmat_{2,2} \Fvec_2
\end{equation}
where $\Uvec$ is the vector of approximate solution values $U_q$ over
the fine nodes in $\tvec_2$.  The MRSDC scheme used in {\tt SMC}
has two operating modes:
\begin{enumerate}
\item The coarse ($\fslow$) and fine ($\ffast$) components
  of~\eqref{eq:mrode} correspond to the advection-diffusion and
  reaction parts of \eqref{eq:cns}, respectively, and the
  forward-Euler sub-stepping scheme is used for both components.
\item The coarse ($\fslow$) and fine ($\ffast$) components
  of~\eqref{eq:mrode} correspond to the reaction and
  advection-diffusion parts of \eqref{eq:cns}, respectively.  The
  advection-diffusion component is treated with a forward-Euler
  sub-stepper while the reaction component is treated with the
  variable-order BDF scheme in VODE \cite{vode}.
\end{enumerate}
In the first mode, the scheme is purely explicit whereas in the second mode
the algorithm uses an implicit treatment of reactions.  Which version is most
appropriate for a given simulation depends on the relative stiffness of chemistry and
fluid dynamics.  When the chemical time scales are not significantly disparate
from the fluid mechanical time scale, the fully explicit version will be the most
efficient. When the chemistry is extremely stiff, the second mode is preferred.

The MRSDC update equation used in {\tt SMC} for mode 1 over the nodes
in $\tvec_2$ is given by
\begin{equation}
  \label{eq:mrupdate}
  \begin{aligned}
  U^{k+1}_{q+1} = U^{k+1}_q
    &+ \Delta t_q
      \bigl[ \fslow(U^{k+1}_{p}, t_{p}) - \fslow(U^k_{p}, t_{p}) \bigr] \\
    &+ \Delta t_q
      \bigl[ \ffast(U^{k+1}_{q}, t_{q}) - \ffast(U^k_{q}, t_{q}) \bigr]
      + \Delta t\, S_{2,2}^{k,q} + \Delta t\, S_{2,1}^{k,q}
  \end{aligned}
\end{equation}
where $S_{2,j}^{k,q}$ is the $q^{\rm th}$ row of $\Smat_{2,j}
\Fvec_j$, and the $p^{\rm th}$ node of the fine component corresponds
to the closest coarse node to the left of the fine node $q$.  That is,
the coarse component $\fslow$ is held constant, and equal to its value
at the left coarse node $p = \lfloor q / (M_2 / M_1) \rfloor$, during
the fine updating steps.
The update equation for mode 2 is
\begin{equation}
  \label{eq:mrupdate2}
  \begin{aligned}
  U^{k+1}_{q+1} = U^{k+1}_q
    &+ \Delta t_q
      \bigl[ \fslow(U^{k+1}_{p+1}, t_{p}) - \fslow(U^k_{p+1}, t_{p}) \bigr] \\
    &+ \Delta t_q
      \bigl[ \ffast(U^{k+1}_{q}, t_{q}) - \ffast(U^k_{q}, t_{q}) \bigr]
      + \Delta t\, S_{2,2}^{k,q} + \Delta t\, S_{2,1}^{k,q}
  \end{aligned}
\end{equation}
where $\fslow(U^{k+1}_{p+1}, t_{p}) $ is computed with VODE.

The computational cost of an MRSDC scheme is determined by the number
of nodes $M_\ell+1$ chosen for each multirate component $\ell$ and the
number of SDC iterations $\Niters$ taken per time step.  The number of
resulting function evaluations $\Nfeval^\ell$ per component $\ell$ is
$\Nfeval^\ell = K M_\ell$.  For example, with 3 Gauss-Lobatto ($M=2$)
nodes for the coarse component and 5 Gauss-Lobatto nodes for the fine
component between each of the coarse nodes (so that there are 9 nodes
in $\tvec_2$ and hence $M_2=8$), the number of coarse function
evaluations over $\Niters=4$ iterations is $\Nfeval^1 = 8$, and the
number of fine function evaluations is $\Nfeval^2 = 32$.

The reacting compressible Navier-Stokes equations \eqref{eq:cns} have
two characteristic time-scales: one over which advection and diffusion
processes evolve, and another over which reactions evolve.  The
reaction time-scale is usually shorter than the advection/diffusion
time-scale.  As such, the coarse and fine MRSDC nodes typically
correspond to the advection/diffusion and reaction parts of
\eqref{eq:cns}, respectively.  However, as is the case for the
Dimethyl Ether Jet problem shown in \S\ref{sect:dme}, the reaction
terms in \eqref{eq:cns} can become so stiff that it becomes
impractical to integrate them explicitly.  In these cases, {\tt SMC}
employs the variable-order BDF schemes in VODE \cite{vode} to
integrate the reaction terms implicitly.  This effectively reverses
the time-step constraint of the system: the reaction terms can be
integrated using VODE (which internally takes many steps) with a
larger time-step than the advection-diffusion terms.  As such, the
roles of the coarse and fine MRSDC nodes can be reversed so that the
coarse and fine MRSDC nodes correspond to the reaction and
advection-diffusion terms of \eqref{eq:cns}, respectively, thereby
amortizing the cost of the implicit solves by using fewer SDC nodes to
integrate the stiff chemistry.

%% file: parallel.tex
{\tt SMC} is implemented within the BoxLib software framework with a
hybrid MPI-OpenMP approach for parallelization \cite{BoxLib}.  In
BoxLib, the computational domain is divided into boxes.  Each MPI
process owns one or more boxes via decomposition of the spatial domain
and it can spawn multiple OpenMP threads.  For parallelization with
OpenMP, the traditional approach is a fine-grained approach in which
individual loops are threaded using OpenMP PARALLEL DO directives
(with COLLAPSE clause if the number of threads is relatively large).
This approach works well for tasks like computing transport
coefficients and chemistry production rates because there is only one
PARALLEL DO region for each box.  However, this approach to
fine-grained parallelism does not work well for computing diffusive
and hyperbolic terms.  The reason for this is that many loops are used
in computing terms like $\nabla \cdot \mathcal{F}_k$.  This results in
many PARALLEL DO regions with serial parts in between that, in turn,
can incur significant overheads and limit parallel efficiency.  To
overcome these issues in the fine-grained parallelism approach, a
coarser-grained approach is taken for computing diffusive and
hyperbolic terms.  In this approach, OpenMP PARALLEL DO directives are
not used.  Instead, a SPMD (single program, multiple data) style is
used.  Each thread works on its own domain.
Listing~\ref{lst:diffterm} shows a code snippet illustrating the
coarse-grained parallelism approach.  Here {\tt q}, {\tt D} and {\tt
  dUdt} are four-dimensional Fortran arrays for primitive variables,
transport coefficients, and rate of change of conserved variables,
respectively.  These shared arrays contain data for the entire box.
The {\tt get\_threadbox} subroutine computes the bounds of sub-boxes
owned by each thread and stores them in private arrays {\tt lo} and
{\tt hi}.  The original box (without ghost cells) is entirely covered
by the sub-boxes without overlap.  Typically, the domain decomposition
among threads is performed in $y$ and $z$-dimensions, but not in
$x$-dimension, because the $x$-dimension of our data arrays is
contiguous in memory.  The {\tt diffterm} subroutine is called by each
thread with private {\tt lo} and {\tt hi}, and it has local variables
(including arrays) that are automatically private to each thread.
Note that there are no OpenMP statements in the {\tt diffterm}
subroutine.  This approach reduces synchronization points and thread
overhead, resulting in improved thread scaling as demonstrated in the
next section.

\begin{lstlisting}[caption={Fortran code snippet illustrating the
coarse-grained parallelism approach.}, label=lst:diffterm]
subroutine compute_dUdt()  ! for a box of nx * ny * nz cells
  ! Only part of the subroutine is shown here
  integer :: lo(3), hi(3)
  !$omp parallel private(lo,hi)
  call get_threadbox(lo,hi)
  call diffterm(lo,hi,q,D,dUdt)
  !$omp end parallel
end subroutine compute_dUdt

subroutine diffterm(lo,hi,q,D,dUdt)
  ! Only part of the subroutine is shown here
  integer, intent(in) :: lo(3), hi(3)
  ! q and D have four ghost cells in each side
  real, intent(in)    ::    q(-3:nx+4,-3:ny+4,-3:nz+4,nq)
  real, intent(in)    ::    D(-3:nx+4,-3:ny+4,-3:nz+4,nD)
  real, intent(inout) :: dUdt( 1:nx  , 1:ny  , 1:nz  ,nU)
  ! tmp is a local array
  real :: tmp(lo(1)-4:hi(1)+4,lo(2)-4:hi(2)+4,lo(3)-4:hi(3)+4)
  ! a nested loop
  ! In the fine-grained approach, this line would be !$omp do.
  do k    =lo(3)-4,hi(3)+4
    do j  =lo(2)-4,hi(2)+4
      do i=lo(1)-4,hi(1)+4
        tmp(i,j,k) = ......
      end do
    end do
  end do
  ! another nested loop
  ! In the fine-grained approach, this line would be !$omp do.
  do k    =lo(3),hi(3)
    do j  =lo(2),hi(2)
      do i=lo(1),hi(1)
        dUdt(i,j,k,1) = ......
      end do
    end do
  end do
end subroutine diffterm
\end{lstlisting}

%% file: results.tex
In this section we present numerical tests to demonstrate the theoretical rate of
convergence and performance of {\tt SMC}.  We also compare
simulation results obtained using {\tt SMC} with the results from the
low Mach number
{\tt LMC} code \cite{DayBell2000,BellDay2001a,nonaka2012deferred}.
Numerical simulations presented here use EGLIB
\cite{Ern1994,Ern1995} to compute transport coefficients.

\subsection{Convergence}
\label{sect:convergence}

\subsubsection{Space/time Convergence}

In this section we present numerical tests to demonstrate the theoretical rate of
convergence of {\tt SMC} for the single-rate SDC integrator.

The convergence tests presented here are set up as follows.  The
simulations use a 9-species $\mathrm{H_2}/\mathrm{O_2}$ reaction
mechanism \cite{LiDryer}.  The computational domain is $L_x = L_y =
L_z = (-0.0005\,\mathrm{m}, 0.0005\,\mathrm{m})$ with periodic
boundaries in all three dimensions.  The initial pressure, temperature
and velocity of the gas are set to
\begin{align}
  p &= p_0 \left[1 + 0.1 \exp{\left(-\frac{r^2}{r_0^2}\right)}\right], \\
  T &= T_0 + T_1 \exp{\left(-\frac{r^2}{r_0^2}\right)},\\
  v_x &= v_0 \sin{\left(\frac{2\pi}{L_x} x\right)}
            \cos{\left(\frac{2\pi}{L_y} y\right)}
            \cos{\left(\frac{2\pi}{L_z} z\right)}, \\
  v_y &= -v_0 \cos{\left(\frac{2\pi}{L_x} x\right)}
             \sin{\left(\frac{2\pi}{L_y} y\right)}
             \cos{\left(\frac{2\pi}{L_z} z\right)}, \\
  v_z &= 0,
\end{align}
where $p_0 = 1\,\mathrm{atm}$, $T_0=1100\,\mathrm{K}$, $T_1 =
400\,\mathrm{K}$, $v_0 = 3\,\mathrm{m}\,\mathrm{s}^{-1}$, $r_0 =
0.0001\,\mathrm{m}$, and $r = \sqrt{x^2+y^2+z^2}$.  The mole fractions
of species are initially set to zero except that
\begin{align}
X(\mathrm{H}_2) &= 0.1 + 0.025 \exp{\left(-\frac{r^2}{r_0^2}\right)}, \\
X(\mathrm{O}_2) &= 0.25 + 0.050 \exp{\left(-\frac{r^2}{r_0^2}\right)}, \\
X(\mathrm{N}_2) &= 1 - X(\mathrm{H}_2) - X(\mathrm{O}_2).
\end{align}
To demonstrate convergence, we performed a series of runs with increasing
spatial and temporal resolutions.  A single-rate SDC scheme with 5
Gauss-Lobatto nodes and 8 SDC iterations is used for time stepping.
The simulations are stopped at $t = 8 \times 10^{-7}\,\mathrm{s}$.
When the spatial resolution changes by a factor of 2 from one run to
another, we also change the time step by a factor of 2.  Therefore,
the convergence rate in this study is expected to be $8^{\rm th}$
order.  Table~\ref{tab:srsdc_cons} shows the $L_\infty$ and $L_2$-norm
errors and convergence rates for $\rho$, $T$, $\vec{u}$,
$Y(\mathrm{H_2})$, $Y(\mathrm{O_2})$, $Y(\mathrm{OH})$,
$Y(\mathrm{H_2O})$, and $Y(\mathrm{N_2})$.  Because there is no
analytic solution for this problem, we compute the errors by comparing
the numerical solution using $N^3$ points with the solution using
$(2N)^3$ points.  The $L_2$-norm error is computed as
\begin{equation}
  L_{2}(N) = \sqrt{\frac{1}{N^3} \sum_{i=1}^N \sum_{j=1}^N \sum_{k=1}^N (\phi_{N,ijk} - \phi_{2N,ijk})^2},
\end{equation}
where $\phi_{N,ijk}$ is the solution at point $(i,j,k)$ of the $N^3$
point run and $\phi_{2N,ijk}$ is the solution of the $(2N)^3$ point
run at the same location.  The results of our simulations demonstrate
the high-order convergence rate of our scheme.  The $L_\infty$ and
$L_2$-norm convergence rates for various variables in the simulations
approach $8^{\rm th}$ order, the theoretical order of the scheme.
Relatively low orders are observed at low resolution because the
numerical solutions are not in the regime of asymptotic convergence
yet.

\begin{table}
  \tbl{Errors and convergence rates for three-dimensional simulations
    of a hydrogen flame using single-rate SDC and 8$^{\rm th}$ order
    finite difference stencil.  Density $\rho$, temperature $T$ and
    $\vec{u}$ are in SI units.\\}{
  \begin{tabular}{lcccccc} \toprule
    Variable & No. of Points & $\Delta t\,(10^{-9}\,\mathrm{s})$ & $L_\infty$ Error &
    $L_\infty$ Rate & $L_2$ Error & $L_2$ Rate \\ \midrule
    $\rho$ & $32^3$  & 4   & 2.613E-08 &      & 5.365E-09 &      \\
           & $64^3$  & 2   & 3.354E-10 & 6.28 & 6.655E-11 & 6.33 \\
           & $128^3$ & 1   & 2.053E-12 & 7.35 & 3.438E-13 & 7.60 \\
           & $256^3$ & 0.5 & 8.817E-15 & 7.86 & 1.437E-15 & 7.90 \\
           \\
    $T$    & $32^3$  & 4    & 3.328E-02 &      & 6.722E-03 &      \\
           & $64^3$  & 2    & 4.184E-04 & 6.31 & 8.266E-05 & 6.35 \\
           & $128^3$ & 1    & 2.515E-06 & 7.38 & 4.266E-07 & 7.60 \\
           & $256^3$ & 0.5  & 1.095E-08 & 7.84 & 1.783E-09 & 7.90 \\
           \\
    $u_x$  & $32^3$  & 4   & 6.184E+00 &      & 8.451E-01 &      \\
           & $64^3$  & 2   & 7.368E-02 & 6.39 & 8.859E-03 & 6.58 \\
           & $128^3$ & 1   & 4.554E-04 & 7.34 & 4.559E-05 & 7.60 \\
           & $256^3$ & 0.5 & 2.011E-06 & 7.82 & 1.904E-07 & 7.90 \\
           \\
    $u_y$  & $32^3$  & 4   & 6.677E+00 &      & 8.801E-01 &      \\
           & $64^3$  & 2   & 7.218E-02 & 6.53 & 9.121E-03 & 6.59 \\
           & $128^3$ & 1   & 4.784E-04 & 7.24 & 4.693E-05 & 7.60 \\
           & $256^3$ & 0.5 & 2.041E-06 & 7.87 & 1.961E-07 & 7.90 \\
           \\
    $u_z$  & $32^3$  & 4   & 6.443E+00 &      & 8.632E-01 &      \\
           & $64^3$  & 2   & 7.211E-02 & 6.48 & 8.987E-03 & 6.59 \\
           & $128^3$ & 1   & 4.626E-04 & 7.28 & 4.624E-05 & 7.60 \\
           & $256^3$ & 0.5 & 2.034E-06 & 7.83 & 1.932E-07 & 7.90 \\
           \\
    $Y(\mathrm{H_2})$ & $32^3$  & 4   & 5.507E-08 &      & 1.272E-08 & \\
                      & $64^3$  & 2   & 8.619E-10 & 6.00 & 1.764E-10 & 6.17 \\
                      & $128^3$ & 1   & 5.802E-12 & 7.21 & 9.828E-13 & 7.49 \\
                      & $256^3$ & 0.5 & 2.512E-14 & 7.85 & 4.169E-15 & 7.88 \\
           \\
    $Y(\mathrm{O_2})$ & $32^3$  & 4   & 1.358E-06 &      & 9.166E-08 & \\
                      & $64^3$  & 2   & 1.406E-08 & 6.59 & 7.136E-10 & 7.00 \\
                      & $128^3$ & 1   & 7.111E-11 & 7.63 & 3.468E-12 & 7.68 \\
                      & $256^3$ & 0.5 & 2.982E-13 & 7.90 & 1.433E-14 & 7.92 \\
           \\
    $Y(\mathrm{OH})$ & $32^3$  & 4   & 3.858E-12 &      & 3.791E-14 & \\
                     & $64^3$  & 2   & 2.740E-14 & 7.14 & 2.661E-16 & 7.15 \\
                     & $128^3$ & 1   & 1.288E-16 & 7.73 & 1.245E-18 & 7.74 \\
                     & $256^3$ & 0.5 & 5.278E-19 & 7.93 & 5.093E-21 & 7.93 \\
           \\
    $Y(\mathrm{H_2O})$ & $32^3$  & 4   & 7.994E-12 &      & 7.318E-14 & \\
                       & $64^3$  & 2   & 5.847E-14 & 7.09 & 5.347E-16 & 7.10 \\
                       & $128^3$ & 1   & 2.778E-16 & 7.72 & 2.548E-18 & 7.71 \\
                       & $256^3$ & 0.5 & 1.141E-18 & 7.93 & 1.048E-20 & 7.93 \\
           \\
    $Y(\mathrm{N_2})$ & $32^3$  & 4   & 1.328E-06 &      & 8.762E-08 & \\
                      & $64^3$  & 2   & 1.383E-08 & 6.58 & 6.622E-10 & 7.05 \\
                      & $128^3$ & 1   & 6.999E-11 & 7.63 & 3.160E-12 & 7.71 \\
                      & $256^3$ & 0.5 & 2.943E-13 & 7.89 & 1.300E-14 & 7.92 \\
    \bottomrule \\
  \end{tabular}}
  \label{tab:srsdc_cons}
\end{table}

\subsubsection{Multirate Integration}

In this section we present numerical tests to demonstrate the theoretical
convergence rate of the MRSDC integrator in {\tt SMC}.

Table~\ref{tab:mrsdc_conv} shows the $L_\infty$ and $L_2$-norm errors and
convergence rates for $\rho$, $T$, $u_x$, $Y(\mathrm{H_2})$,
$Y(\mathrm{O_2})$, $Y(\mathrm{OH})$, $Y(\mathrm{H_2O})$, and
$Y(\mathrm{N_2})$ for the same test problem as described in
\S\ref{sect:convergence}, except that $T_0=300\,\mathrm{K}$ and $T_1 =
1100\,\mathrm{K}$.  All runs were performed on a $32^3$ grid,
so that the errors and convergence rates reported here are ODE errors.
Because there is no analytic solution for this problem, a reference
solution was computed using the same spatial grid, but with an $8^{\rm
  th}$-order single-rate SDC integrator and $\Delta t = 1.25 \times 10^{-9}$s.  We
note that the $L_\infty$ and $L_2$-norm convergence rates for all variables
and MRSDC configurations are $4^{\rm th}$-order, as expected since the
coarse MRSDC component uses 3 Gauss-Lobatto nodes.  This verifies that
the MRSDC implementation in {\tt SMC} operates according to its
specifications regardless of how the fine nodes are chosen (see
\S\ref{sect:mrsdc}).

%
%
\begin{table}
  \tbl{Errors and convergence rates for three-dimensional simulations
    of a hydrogen flame using 4$^{\rm th}$ order MRSDC and 8$^{\rm th}$ order
    finite difference stencil.  Density $\rho$, temperature $T$ and
    $\vec{u}$ are in SI units.
    The MRSDC $X$ / $Y$ configurations use $X$ coarse nodes and $Y$ fine nodes between each pair of coarse nodes (type (b) in \S\ref{sect:mrsdc}).
    The MRSDC $X$ / $Y \times R$ configurations use $X$ coarse nodes
    and $Y$ fine nodes repeated $R$ times between each pair of coarse
    nodes (type (c) in \S\ref{sect:mrsdc}).  All SDC nodes are
    Gauss-Lobatto (GL) nodes. Density $\rho$, temperature $T$ and
    $\vec{u}$ are in SI units.\\}{
  \begin{tabular}{llccccc} \toprule
SDC Configuration & Variable & $\Delta t\,(10^{-9}\,\mathrm{s})$ & $L_\infty$ Error &
    $L_\infty$ Rate & $L_2$ Error & $L_2$ Rate \\ \midrule
GL 3 / GL 9  & $\rho$ & 10  & 5.638e-07 &      & 3.958e-08 & \\
             &        & 5   & 3.392e-08 & 4.06 & 2.426e-09 & 4.03 \\
             &        & 2.5 & 2.093e-09 & 4.02 & 1.509e-10 & 4.01 \\
\\
 & $T$ & 10 & 2.748e-04 & & 1.577e-05 & \\
 &     & 5 & 1.661e-05 & 4.05 & 9.504e-07 & 4.05 \\
 &     & 2.5 & 1.029e-06 & 4.01 & 5.877e-08 & 4.02 \\
\\
 & $u_x$ & 10  & 6.928e-04 &      & 2.253e-05 & \\
 &       & 5   & 4.255e-05 & 4.03 & 1.369e-06 & 4.04 \\
 &       & 2.5 & 2.650e-06 & 4.01 & 8.496e-08 & 4.01 \\
\\
 & $Y({\rm H}_2)$ & 10 & 1.960e-09 & & 5.051e-11 & \\
 &  & 5 & 1.160e-10 & 4.08 & 3.059e-12 & 4.05 \\
 &  & 2.5 & 7.108e-12 & 4.03 & 1.896e-13 & 4.01 \\
\\
 & $Y({\rm O}_2)$ & 10 & 9.501e-09 & & 3.340e-10 & \\
 &  & 5 & 5.762e-10 & 4.04 & 2.027e-11 & 4.04 \\
 &  & 2.5 & 3.562e-11 & 4.02 & 1.255e-12 & 4.01 \\
\\
 & $Y({\rm OH})$ & 10 & 5.989e-17 & & 3.531e-19 & \\
 &  & 5 & 3.671e-18 & 4.03 & 2.183e-20 & 4.02 \\
 &  & 2.5 & 2.292e-19 & 4.00 & 1.370e-21 & 3.99 \\
\\
 & $Y({\rm H}_2{\rm O})$ & 10 & 4.309e-18 & & 6.660e-20 & \\
 &  & 5 & 2.519e-19 & 4.10 & 3.928e-21 & 4.08 \\
 &  & 2.5 & 1.527e-20 & 4.04 & 2.394e-22 & 4.04 \\
\\
 & $Y({\rm N}_2)$ & 10 & 8.743e-09 & & 3.191e-10 & \\
 &  & 5 & 5.258e-10 & 4.06 & 1.938e-11 & 4.04 \\
 &  & 2.5 & 3.238e-11 & 4.02 & 1.200e-12 & 4.01 \\
\\
GL 3 / GL 5 $\times$ 2 & $\rho$ & 10  & 5.638e-07 &      & 3.958e-08 & \\
                       &        & 5   & 3.392e-08 & 4.06 & 2.426e-09 & 4.03 \\
                       &        & 2.5 & 2.093e-09 & 4.02 & 1.509e-10 & 4.01 \\
\\
 & $T$ & 10 & 2.748e-04 & & 1.577e-05 & \\
 &  & 5 & 1.661e-05 & 4.05 & 9.504e-07 & 4.05 \\
 &  & 2.5 & 1.029e-06 & 4.01 & 5.877e-08 & 4.02 \\
\\
 & $u_x$ & 10  & 6.928e-04 &      & 2.253e-05 & \\
 &       & 5   & 4.255e-05 & 4.03 & 1.369e-06 & 4.04 \\
 &       & 2.5 & 2.650e-06 & 4.01 & 8.496e-08 & 4.01 \\
\\
 & $Y({\rm H}_2)$ & 10 & 1.960e-09 & & 5.051e-11 & \\
 &  & 5 & 1.160e-10 & 4.08 & 3.059e-12 & 4.05 \\
 &  & 2.5 & 7.108e-12 & 4.03 & 1.896e-13 & 4.01 \\
\\
 & $Y({\rm O}_2)$ & 10 & 9.501e-09 & & 3.340e-10 & \\
 &  & 5 & 5.762e-10 & 4.04 & 2.027e-11 & 4.04 \\
 &  & 2.5 & 3.562e-11 & 4.02 & 1.255e-12 & 4.01 \\
\\
 & $Y({\rm OH})$ & 10 & 5.989e-17 & & 3.531e-19 & \\
 &  & 5 & 3.671e-18 & 4.03 & 2.183e-20 & 4.02 \\
 &  & 2.5 & 2.292e-19 & 4.00 & 1.370e-21 & 3.99 \\
\\
 & $Y({\rm H}_2{\rm O})$ & 10 & 4.309e-18 & & 6.660e-20 & \\
 &  & 5 & 2.519e-19 & 4.10 & 3.928e-21 & 4.08 \\
 &  & 2.5 & 1.527e-20 & 4.04 & 2.394e-22 & 4.04 \\
\\
 & $Y({\rm N}_2)$ & 10 & 8.743e-09 & & 3.191e-10 & \\
 &  & 5 & 5.258e-10 & 4.06 & 1.938e-11 & 4.04 \\
 &  & 2.5 & 3.238e-11 & 4.02 & 1.200e-12 & 4.01 \\ \bottomrule
  \end{tabular}}
  \label{tab:mrsdc_conv}
\end{table}

\subsection{Performance}
\label{sect:performance}

\subsubsection{Parallel Performance}

In this section we present numerical tests to demonstrate the parallel
efficiency of {\tt SMC}.

Two strong scaling studies were performed: a pure OpenMP study on a
61-core Intel Xeon Phi, and a pure MPI study on the Hopper
supercomputer at the National Energy Research Scientific Computing
Center (NERSC).

The pure OpenMP strong scaling study was carried out on a 61-core 1.1
GHz Intel Xeon Phi coprocessor that supported up to 244 hyperthreads.
We performed a series of pure OpenMP runs, each consisting of a single
box of $128^3$ points, of the test problem in
section~\ref{sect:convergence} using various numbers of threads.  The
simulations are run for 10 time steps.  Table~\ref{tab:strong-scaling}
shows that {\tt SMC} has excellent thread scaling behavior.  An factor of 86
speedup over the 1 thread run is obtained on the 61-core
coprocessor using 240 hyperthreads.

The pure MPI strong scaling study was carried out on the Hopper
supercomputer at NERSC.  Each compute node of Hopper has 4 non-uniform
memory access (NUMA) nodes, with 6 cores on each NUMA node.  We
performed a series of pure MPI runs of the test problem in
section~\ref{sect:convergence} using various numbers of processors and
integration schemes on a grid of $64^3$ points.  The simulations are
run for 10 time steps.  Several observations can be made from the
results listed in Table~\ref{tab:strong-scaling-mpi}.  First, when
using a single processor the wide stencil is slightly faster than the
narrow stencil regardless of the integration scheme used.  This is
consistent with the cost analysis of the stencils done in
\S\ref{sec:cost} which shows that the wide stencil uses fewer FLOPs
per stencil evaluation than the narrow stencil.  The difference in
single-processor run times, however, is not significant which suggests
that both stencils are operating in a regime dominated by memory
access times.  Second, when using multiple MPI processes, both the
wall clock time per step and efficiency of the narrow stencil is
consistently better than the wide stencil regardless of the integrator
being used.  Again, this is consistent with cost analysis done in
\S\ref{sec:cost} which shows that the communication cost associated
with the narrow stencil is less than the wide stencil.  The difference
between the stencils when using multiple MPI processes is observable
but not dramatic.  Regarding parallel efficiency, we note that as more
processors are used the amount of work performed grows due to the
presence of ghost cells where, for example, diffusion coefficients
need to be computed.  As such, we cannot expect perfect parallel
efficiency.  In the extreme case were the problem is spread across 512
processors, each processor treats an $8^3$ box, which grows to a
$16^3$ box with the addition of ghost cells.  Finally, note that the
parallel efficiency of the MRSDC scheme is higher than the SRSDC
scheme regardless of the stencil being used since, for this test, the
MRSDC scheme performs the same number of advection/diffusion
evaluations (and hence stencil evaluations and extra ghost cell work)
but performs significantly more evaluations of the chemistry terms --
MRSDC has a higher arithmetic intensity with respect to MPI
communications than SRSDC.

\begin{table}
\tbl{Strong scaling behavior of pure OpenMP runs on a 61-core Intel
  Xeon Phi coprocessor.  Average wall clock time per time step and speedup
  are shown for runs using various number of threads.  Hyperthreading
  is used for the 128 and 240-thread runs.\\}{
\begin{tabular}{lccc} \toprule
 No. of threads & Wall time per step (s) & Speedup \\
 \midrule
     1 & 1223.58 &      \\
     2 &  626.13 & 2.0 \\
     4 &  328.97 & 3.7 \\
     8 &  159.88 & 7.7 \\
    16 &   79.73 & 16 \\
    32 &   40.12 & 31 \\
    60 &   24.55 & 50 \\
   128 &   18.19 & 67 \\
   240 &   14.28 & 86 \\
 \bottomrule \\
  \end{tabular}}
  \label{tab:strong-scaling}
\end{table}

\begin{table}
\tbl{Strong scaling behavior of pure MPI runs with $64^3$ grid points on the Hopper supercomputer at NERSC.  Average wall clock time per time step, speedup, and efficiency
  are shown for runs using various number of processors and schemes.
  The SRSDC scheme uses single-rate SDC with 3 nodes, with narrow and wide stencils.
  The MRSDC scheme uses multi-rate SDC with 3 fine nodes, and 5 fine nodes repeated 2 times between each coarse node (type (c) in \S\ref{sect:mrsdc}).  All SDC nodes are Gauss-Lobatto (GL) nodes.
\\}{
\begin{tabular}{llcccc} \toprule
Method/Stencil & No. of MPI processes & Wall time per step (s) & Speedup & Efficiency \\
 \midrule
 SRSDC/Narrow & 1   & 243.35 &        &      \\
              & 8   &  42.92 & 5.67   & 71\% \\
              & 64  &   7.91 & 30.84  & 48\% \\
              & 512 &   2.14 & 133.95 & 22\% \\
 SRSDC/Wide   & 1   & 231.60 &        &      \\
              & 8   &  43.21 & 5.36   & 67\% \\
              & 64  &   8.04 & 28.80  & 45\% \\
              & 512 &   2.27 & 102.09 & 20\% \\
 MRSDC/Narrow & 1   & 695.57 &        &      \\
              & 8   & 101.86 & 6.82   & 85\% \\
              & 64  &  15.25 & 45.69  & 71\% \\
              & 512 &   3.15 & 220.55 & 43\% \\
 MRSDC/Wide   & 1   & 685.87 &        &      \\
              & 8   & 102.45 & 6.69   & 84\% \\
              & 64  &  15.75 & 43.55  & 68\% \\
              & 512 &   3.28 & 208.85 & 41\% \\
 \bottomrule\\
  \end{tabular}}
  \label{tab:strong-scaling-mpi}
\end{table}

We also performed a weak scaling study using the test problem in
section~\ref{sect:convergence} on the Hopper supercomputer at the
NERSC.  In this weak scaling study, a hybrid MPI/OpenMP approach with
6 threads per MPI process is used.  Each MPI process works on two
$64^3$ boxes.  The simulations are run for 10 time steps.
Table~\ref{tab:weak-scaling} shows the average wall clock time per
time step for a series of runs on various numbers of cores.  We use the
6-core run as baseline and define parallel efficiency for an $n$-core
run as $T_6/T_n$, where $T_6$ and $T_n$ are the average wall clock
time per time step for the 6 and $n$-core runs, respectively.
Excellent scaling to about 100 thousand cores is observed.

\begin{table}
\tbl{Weak scaling behavior of hybrid MPI/OpenMP runs on Hopper at NERSC.  Average
  wall clock time per time step and parallel efficiency are shown for
  runs on various number of cores.  Each MPI process spawns 6 OpenMP
  threads.\\}{
\begin{tabular}{lccc} \toprule
 No. of cores & No. of MPI processes & Wall time per step (s) & Efficiency \\
 \midrule
     6 &     1 & 13.72 & \\
    24 &     4 & 14.32 & 96\% \\
   192 &    32 & 14.46 & 95\% \\
  1536 &   256 & 14.47 & 95\% \\
 12288 &  2048 & 15.02 & 91\% \\
 98304 & 16384 & 15.32 & 90\% \\
 \bottomrule \\
  \end{tabular}}
  \label{tab:weak-scaling}
\end{table}

\subsubsection{Multirate Integration}

Table~\ref{tab:mrsdc_runtime} compares run time and function
evaluation counts for premixed methane flame using various
configurations of SDC; MRSDC; and a six-stage, fourth-order
Runge-Kutta integrator.  The configuration is a cubic box ($L_x = L_y
= L_z = 0.008$m) with periodic boundary conditions in the $x$ and $y$
directions, and inlet/outlet boundary conditions in the $z$ direction.
The GRI-MECH 3.0 chemistry network (53 species, 325 reactions)
\cite{SmithGolden} was used so that the reaction terms in
\eqref{eq:cns} were quite stiff and an MRSDC integrator with the
reaction terms on the fine nodes was appropriate.  A reference
solution was computed using a single-rate SDC integrator with a time
step $\Delta t$ of $1 \times 10^{-9}$s to a final time of $4.2 \times
10^{-7}$s.  All simulations were run on a $32^3$ point grid.  For each
integrator, the test was run using successively larger time steps
until the solution no longer matched the reference solution.  The
largest time step that matched the reference solution (to within a
relative $L_\infty$ error of less than $1 \times 10^{-5}$) is
reported.  All SDC integrators use 3 Gauss-Lobatto nodes for the
coarse component and hence all integrators, including the Runge-Kutta
integrators, are formally $4^{\rm th}$ order accurate.  We note that
all of the multi-rate integrators are able to use larger time steps
than the SRSDC and RK integrators due to the stiffness of the reaction
terms: the multi-rate integrators use more nodes to integrate the
reaction terms and hence the effective time-steps used in the update
equation \eqref{eq:mrupdate} are smaller.  This means that, in all
cases, the multi-rate integrators evaluate the advection/diffusion
terms less frequently than the single-rate integrator.  However, in
all but one case, the multi-rate integrators also evaluate the
reaction terms more frequently than the SRSDC and RK integrators.  The
trade-off between the frequency of evaluating the advection/diffusion
vs reaction terms is highlighted by comparing the MRSDC 3 GL / 9 GL
and MRSDC 3 GL / 13 GL runs: the run with 13 fine nodes can use a
slightly larger time-step at the cost of evaluating the reaction terms
more frequently, resulting in only a very slight runtime advantage
over the run using 9 fine nodes and a smaller time-step.  Ultimately
the largest runtime advantage is realized by choosing a multi-rate
configuration that allows a large time-step to be taken without
requiring too many extra calls to the reaction network.  For this test
the MRSDC 3 / $5 \times 2$ configuration achieves this goal and is
able to run to completion approximately 6.0 times faster than the
single-rate SDC integrator and 4.5 times faster than the RK
integrator.

%
%

\begin{table}
  \tbl{%
    Comparison of runtime and function evaluation count for various configurations of SDC and MRSDC.
    Function evaluation counts are reported as: no. of advection/diffusion evaluations / no. of chemistry evaluations.
    All schemes use the narrow stencil, except where noted.
    The SRSDC $X$ configuration uses single-rate SDC with $X$ nodes.
    The RK scheme is fourth order and use six stages.
    The MRSDC $X$ / $Y$ configurations use MRSDC with $X$ coarse nodes and $Y$ fine nodes between each pair of coarse nodes (type (b) in \S\ref{sect:mrsdc}).
    The MRSDC $X$ / $Y \times R$ configurations use MRSDC with $X$ coarse nodes and $Y$ fine nodes repeated $R$ times between each pair of coarse nodes (type (c) in \S\ref{sect:mrsdc}).  All SDC nodes are Gauss-Lobatto (GL) nodes.\\}
  {\begin{tabular}{lcrr} \toprule
    Scheme & $\Delta t$ ($10^{-9}$s) & Runtime (s) & Evaluations (AD/R) \\ \midrule
    SRSDC, 3 GL & 6 & 130.6 & 561/561 \\
    RK, 6/4 & 6 & 96.5 & 420/420 \\
    RK, 6/4 wide & 6 & 95.4 & 420/420 \\
    MRSDC, 3 GL / 9 GL type (b) & 30 & 44.1 & 113/897 \\
    MRSDC, 3 GL / 13 GL type (b)& 42 & 39.7 & 81/961 \\
    MRSDC, 3 GL / $3 \times 8$ GL type (c) & 60 & 30.8 & 57/897 \\
    MRSDC, 3 GL / $5 \times 2$ GL type (c) & 60 & 21.7 & 57/449 \\
    \bottomrule \\
  \end{tabular}}
  \label{tab:mrsdc_runtime}
\end{table}

\subsection{Dimethyl Ether Jet}
\label{sect:dme}

We present here a two-dimensional simulation of jet using a 39-species
dimethyl ether (DME) chemistry mechanism \cite{DME39}.
A 2D Cartesian domain, $-0.00114\,\mathrm{m} < x <
0.00114\,\mathrm{m}$ and $0 < y < 0.00228\,\mathrm{m}$, is used
for this test.  Pressure is
initially set to $40\,\mathrm{atm}$ everywhere.  The initial
temperature, velocity and mole fractions of species are set to
\begin{align}
T_0 &= \eta T_{\mathrm{jet}} + (1-\eta) T_{\mathrm{air}}, \\
v_{0x} &= 0, \\
v_{0y} &= \eta v_{\mathrm{jet}} + (1-\eta) v_{\mathrm{air}}, \\
X_0 &= \eta X_{\mathrm{jet}} + (1-\eta) X_{\mathrm{air}},
\end{align}
where $T_{\mathrm{jet}} = 400\,\mathrm{K}$, $T_{\mathrm{air}} =
1525\,\mathrm{K}$, $v_{\mathrm{jet}} =
51.2\,\mathrm{m}\,\mathrm{s}^{-1}$, and $v_{\mathrm{air}} =
5.12\,\mathrm{m}\,\mathrm{s}^{-1}$. The $\eta$ variable is given by
\begin{equation}
  \eta = \frac{1}{2} (\tanh{\frac{x + x_0}{\sigma}} - \tanh{\frac{x - x_0}{\sigma}}),
\end{equation}
where $x_0 = 5.69 \times 10^{-5}\,\mathrm{m}$ and $\sigma = 0.5 x_0$.  The mole
fractions of species for the ``jet'' and ``air'' states are set to
zero except that
\begin{align}
  X_{\mathrm{jet}}(\mathrm{CH_3OCH_3}) &= 0.2, \\
  X_{\mathrm{jet}}(\mathrm{N_2}) &= 0.8, \\
  X_{\mathrm{air}}(\mathrm{O_2}) &= 0.21, \\
  X_{\mathrm{air}}(\mathrm{N_2}) &= 0.79.
\end{align}
An inflow boundary is used at the lower $y$-boundary, whereas outflow
boundary conditions are applied at the other three
boundaries. Sinusoidal variation is added to the inflow velocity as
follows,
\begin{equation}
  v_y(x,t) = v_{0y} + \tilde{v} \eta \sin\left(\frac{2 \pi}{L_x}
    x\right) \sin\left(\frac{2\pi}{L_t} t\right),
\end{equation}
where $\tilde{v} = 10\,\mathrm{m}\,\mathrm{s}^{-1}$, $L_x =
0.00228\,\mathrm{m}$ is the length of the domain in $x$-direction, and
$L_t = 10^{-5}\,\mathrm{s}$ is the period of the variation.
With these parameters, the resulting flow problem is in the low
Mach number regime to facilitate cross-validation with an established
low Mach number solver \cite{DayBell2000}

\begin{figure}
  \centering
  \includegraphics[width=\textwidth]{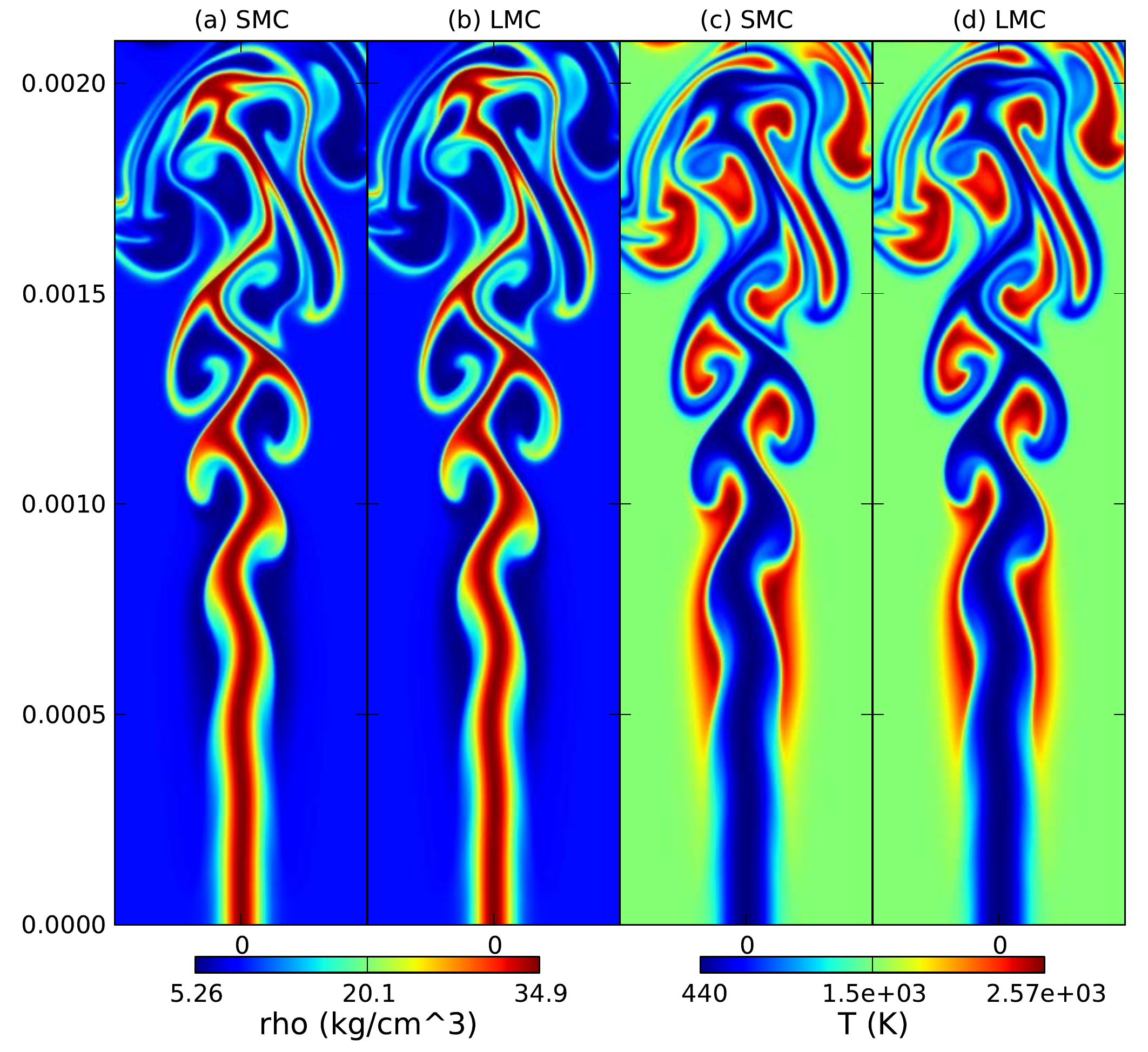}
  \caption{Density and temperature at $6 \times 10^{-5}\,\mathrm{s}$
    for the DME jet test problem. We show (a) density from {\tt SMC},
    (b) density from {\tt LMC}, (c) temperature from {\tt SMC}, and
    (d) temperature from {\tt LMC}.  Only part of the domain
    immediately surrounding the jet is shown here.}
  \label{fig:DME_rhoT}
\end{figure}

\begin{figure}
  \centering
  \includegraphics[width=\textwidth]{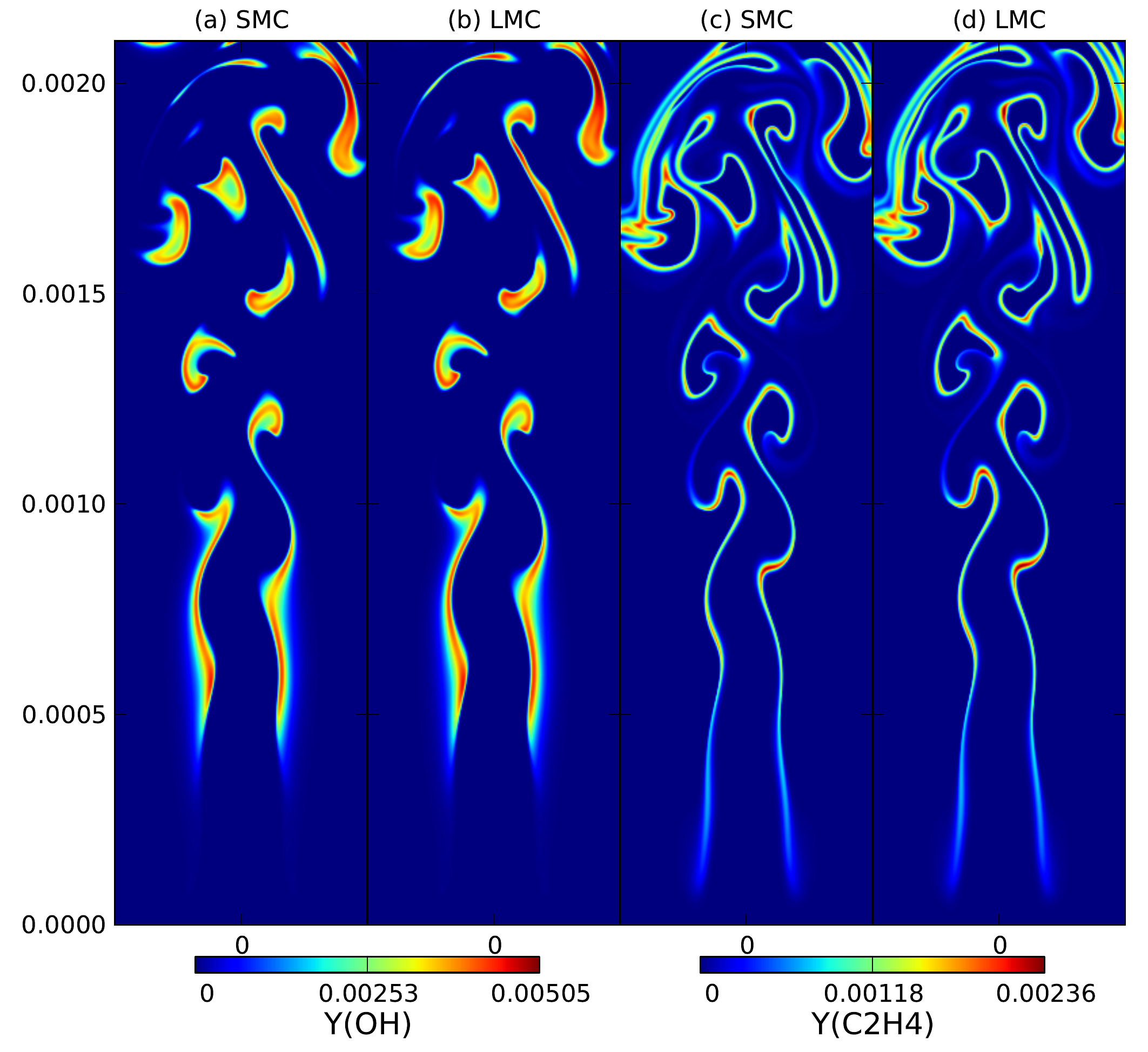}
  \caption{Mass fraction of $\mathrm{OH}$ and $\mathrm{C_2H_4}$ at $6
    \times 10^{-5}\,\mathrm{s}$ for the DME jet test problem. We show
    (a) $Y(\mathrm{OH})$ from {\tt SMC}, (b) $Y(\mathrm{OH})$ from
    {\tt LMC}, (c) $Y(\mathrm{C_2H_4})$ from {\tt SMC}, and (d)
    $Y(\mathrm{C_2H_4})$ from {\tt LMC}.  Only part of the domain
    immediately surrounding the jet is shown here.}
  \label{fig:DME_OH-C2H4}
\end{figure}

The 39-species DME mechanism used in this test is extremely stiff so
that it becomes impractical to evolve the reaction system explicitly.
Hence, we use the variable-order BDF scheme in VODE \cite{vode} for
computing species production rates.  In this test, we use MRSDC, and
put the reaction terms $\rho \dot{\omega}_k$ in (\ref{eq:NS:rY}) and
advection-diffusion terms on the coarse and fine nodes, respectively.
Due to the use of an implicit solver, we can evaluate the reaction
terms on the coarse node with a large time step without suffering from
instability.  As for advection-diffusion terms, the time step is
limited by the acoustic time scale.  In this test, we use 5 coarse
Gauss-Lobatto nodes for the reaction terms, whereas for the
advection-diffusion terms we use 3 fine Gauss-Lobatto nodes repeated
twice between each pair of coarse nodes.  The time step is chosen to
have a Courant-Friedrichs-Lewy (CFL) number of 2.  The effective CFL
numbers for a substep of advection-diffusion (i.e., the interval
spanned by 3 fine Gauss-Lobatto nodes) are about 0.17 and 0.33. (Note
that the 5 coarse Gauss-Lobatto nodes are not evenly distributed.)  We
perform 4 SDC iterations per time step.  The simulation is run with
$2048^2$ uniform points corresponding to a cell size of $\sim 1.1
\times 10^{-6}\,\mathrm{m}$ up to time $6 \times 10^{-5}\,\mathrm{s}$.
The numerical results are shown in figures~\ref{fig:DME_rhoT} \&
\ref{fig:DME_OH-C2H4}.  As a validation of {\tt SMC}, we have also
presented the numerical results obtained using the {\tt LMC} code that
is based on a low-Mach number formulation
\cite{DayBell2000,BellDay2001a,nonaka2012deferred} that approximates
the compressible flow equation to $O(M^2)$.  The results of the two
codes are nearly indistinguishable in spite of the difference between
the two formulations.  This is because the maximum Mach number $M$ in
this problem is only $0.14$ so the low Mach number model is
applicable.  We note that this simulation is essentially infeasible
with a fully explicit solver because of the stiffness of the chemical
mechanism.

%% file: stencil.tex
The truncation error (numerical minus exact) of the $8^{\rm th}$ order
narrow stencil discretization of $\partial(a {\partial u}/{\partial
  x})/{\partial x}$ (\S~\ref{sect:stencil}) is given by
\begin{align}
e_i = {} & -\frac{1}{3150} a \frac{\partial^{10} u}{\partial x^{10}}
 -\frac{1}{630} \frac{\partial^{} a}{\partial x^{}} \frac{\partial^{9} u}{\partial x^{9}}
 -\frac{1}{35} \frac{\partial^{2} a}{\partial x^{2}} \frac{\partial^{8} u}{\partial x^{8}}
 +(-\frac{113}{840} - 4 m_{48})\frac{\partial^{3} a}{\partial x^{3}} \frac{\partial^{7} u}{\partial x^{7}}
\nonumber \\
{} & + (-\frac{487}{1680} - 14 m_{48}) \frac{\partial^{4} a}{\partial x^{4}} \frac{\partial^{6} u}{\partial x^{6}}
 +(-\frac{4513}{12600} + \frac{2}{5} m_{47} - \frac{92}{5} m_{48})\frac{\partial^{5} a}{\partial x^{5}} \frac{\partial^{5} u}{\partial x^{5}}
\nonumber \\
 {} & + (-\frac{2777}{10080} + m_{47} - 11 m_{48}) \frac{\partial^{6} a}{\partial x^{6}} \frac{\partial^{4} u}{\partial x^{4}}
 +(-\frac{3181}{25200} + \frac{4}{5} m_{47} - \frac{14}{5} m_{48}) \frac{\partial^{7} a}{\partial x^{7}} \frac{\partial^{3} u}{\partial x^{3}}
\nonumber \\
 {} & + (-\frac{1403}{50400} + \frac{1}{5}m_{47} - \frac{1}{5}m_{48}) \frac{\partial^{8} a}{\partial x^{8}} \frac{\partial^{2} u}{\partial x^{2}}
- \frac{1}{630} \frac{\partial^{9} a}{\partial x^{9}} \frac{\partial^{} u}{\partial x^{}}
\nonumber \\
        {} & + O(\Delta x^8), \label{eq:truncerror}
\end{align}
where $a$ and all derivatives of $a$ and $u$ are evaluated at $x =
x_i$.  An upper bound of the truncation error can be calculated by
summing the absolute value of each term in \eqref{eq:truncerror}.
Assuming that all $(\partial^na/\partial x^n)(\partial^{10-n}u/\partial
x^{10-n})$ are of similar order, we can minimize the upper bound by
setting $m_{47} = 3557/44100$ and $m_{48} = -2083/117600$.  We do not
expect this particular choice will be the best in general.

We now present numerical results for the two test problems in
\cite{Pantano2010} using the narrow stencil.  The test problems solve the two-dimensional
heat equation
\begin{equation}
  \frac{\partial u}{\partial t} = \frac{\partial}{\partial x} \left(
    a(x,y) \frac{\partial u}{\partial x} \right) +
      \frac{\partial}{\partial y} \left( b(x,y) \frac{\partial
          u}{\partial y} \right) + g(x,y,t),
\end{equation}
on the periodic spatial domain $[0, 2\pi] \times [0, 2 \pi]$.
The initial condition of is given by
\begin{equation}
  u(x,y,t=0) = \sin{(x)} \sin{(y)}.
\end{equation}
In the first test problem (T1), the variable coefficients and source
term are given by
\begin{align}
a(x,y) = { } & b(x,y) = 1 + \epsilon \cos{(x)} \cos{(y)}, \\
g(x,y,t) = { } & (1 + 4 \epsilon \cos{(x)} \cos{(y)}) u_{\mathrm{ex}}(x,y,t),
\end{align}
where $\epsilon = 0.1$, and the exact solution is
\begin{equation}
  u_{\mathrm{ex}}(x,y,t) = e^{-t} \sin{(x)} \sin{(y)}.\label{eq:exact1}
\end{equation}
In the second test problem (T2), the variable coefficients and source
term are given by
\begin{align}
a(x,y) = { } & 1 + \epsilon \cos{(2x)} \sin{(2y+\frac{\pi}{3})}, \\
b(x,y) = { } & 1 + \epsilon \cos{(2x+\frac{\pi}{3})} \sin{(2y)}, \\
g(x,y,t) = { } & 0,
\end{align}
where $\epsilon = 0.9$.  There is no exact analytic solution for the second test.

We have performed a series of calculations with various resolutions.
A $4^{\rm th}$ order SDC scheme with 3 Gauss-Lobatto nodes and 4 SDC
iterations is used with time step $\Delta t = 0.4 ((\Delta x^{-2} +
\Delta y^{-2}) \times \max(a(x,y),b(x,y)))^{-1}$.  Three sets of the
two free parameters in the narrow stencil are selected in these runs:
(1) $m_{47} = 3557/44100$ and $m_{48} = -2083/117600$, hereafter
denoted by ``SMC''; (1) $m_{47} = 0$ and $m_{48} = 0$, hereafter
denoted by ``ZERO''; and (3) $m_{47} = 1059283/13608000$ and $m_{48} =
-856481/40824000$, hereafter denoted by ``OPTIMAL''.  The SMC stencil
minimizes the upper bound of the truncation
error.  
The OPTIMAL stencil is designed to minimize the 2-norm of truncation
error for the test problems given the initial
condition. Tables~\ref{tab:stenciltest1} \& \ref{tab:stenciltest2}
show the $L_\infty$ and $L_2$-norm errors and convergence rates for
the three stencils.  For the second test problem, the errors are
computed by comparing the solutions to a high-resolution run using
$640^2$ cells.  Note that the designed order of convergence is
observed in the numerical experiments.  In particular, the convergence
rate is independent of the choice of the two free parameters, as
expected.  The SMC stencil, which is the default stencil used in {\tt
  SMC}, performs well in both tests.  In contrast, the OPTIMAL
stencil, which is optimized specifically for these two test problems,
produces the smallest errors (especially for the first test problem).
Finally, note that the accuracy of the stencil is essentially
insensitive to the two free parameters for the second, more
challenging, test problem.

\begin{table}
\tbl{Errors and convergence rates for the two-dimensional stencil test
  problem 1 at $t = 1$.  Quadruple precision is used in this test to
  minimize roundoff errors in the floating-point computation.\\}
{
  \begin{tabular}{llcccc} \toprule
    Stencil & No. of Points & $L_\infty$ Error & $L_\infty$ Rate & $L_2$ Error & $L_2$ Rate \\
    \midrule
     SMC &  $20^2$ & 4.167E-08 &      & 1.984E-08 &      \\
       &  $40^2$ & 1.721E-10 & 7.92 & 8.005E-11 & 7.95 \\
       &  $80^2$ & 6.827E-13 & 7.98 & 3.152E-13 & 7.99 \\
       & $160^2$ & 2.676E-15 & 8.00 & 1.234E-15 & 8.00 \\
       \\
     ZERO &  $20^2$ & 2.642E-07 &      & 1.394E-07 &      \\
       &  $40^2$ & 1.223E-09 & 7.76 & 5.911E-10 & 7.88 \\
       &  $80^2$ & 4.877E-12 & 7.97 & 2.357E-12 & 7.97 \\
       & $160^2$ & 1.915E-14 & 7.99 & 9.252E-15 & 7.99 \\
       \\
     OPTIMAL &  $20^2$ & 1.444E-08 &      & 7.320E-09 &      \\
       &  $40^2$ & 5.757E-11 & 7.97 & 2.880E-11 & 7.99 \\
       &  $80^2$ & 2.261E-13 & 7.99 & 1.130E-13 & 7.99 \\
       & $160^2$ & 8.842E-16 & 8.00 & 4.422E-16 & 8.00 \\
    \bottomrule \\
  \end{tabular}
}
\label{tab:stenciltest1}
\end{table}

\begin{table}
\tbl{Errors and convergence rates for the two-dimensional stencil test
  problem 2 at $t = 1$.\\}
{
  \begin{tabular}{llcccc} \toprule
    Stencil & No. of Points & $L_\infty$ Error & $L_\infty$ Rate & $L_2$ Error & $L_2$ Rate \\
    \midrule
     SMC &  $40^2$ & 2.279E-04 &      & 4.270E-05 &  \\
       &  $80^2$ & 4.261E-06 & 5.74 & 6.151E-07 & 6.12 \\
       & $160^2$ & 3.351E-08 & 6.99 & 4.021E-09 & 7.26 \\
       & $320^2$ & 1.639E-10 & 7.68 & 1.835E-11 & 7.78 \\
       \\
     ZERO &  $40^2$ & 2.611E-04 &      & 4.913E-05 &  \\
       &  $80^2$ & 4.744E-06 & 5.78 & 6.918E-07 & 6.15 \\
       & $160^2$ & 3.682E-08 & 7.01 & 4.481E-09 & 7.27 \\
       & $320^2$ & 1.797E-10 & 7.68 & 2.039E-11 & 7.78 \\
       \\
     OPTIMAL &  $40^2$ & 2.221E-04 &      & 4.160E-05 &  \\
       &  $80^2$ & 4.174E-06 & 5.73 & 6.014E-07 & 6.11 \\
       & $160^2$ & 3.291E-08 & 6.99 & 3.938E-09 & 7.25 \\
       & $320^2$ & 1.610E-10 & 7.68 & 1.798E-11 & 7.77 \\
    \bottomrule \\
  \end{tabular}
}
\label{tab:stenciltest2}
\end{table}